\documentclass[11pt,twoside]{amsart}
 
\usepackage{amsfonts}
\usepackage{epsfig,graphics, graphicx}
\usepackage{amssymb}
\usepackage{tikz, tikz-cd}
\usepackage{hyperref, enumerate, indentfirst}
\hypersetup{
  linktocpage,
    colorlinks,
    citecolor=black,
    filecolor=black,
    linkcolor=blue,
    urlcolor=blue
    }

\newtheorem{theorem}{Theorem}[section]

\newtheorem{proposition}[theorem]{Proposition}
       
\newtheorem{definition}[theorem]{Definition}

\newtheorem{lemma}[theorem]{Lemma}

\newtheorem{remark}[theorem]{Remark}

\newtheorem{corollary}[theorem]{Corollary}

\newcommand{\Aff}{\operatorname{Aff }}
\newcommand{\SL}{\operatorname{SL}}

\newcommand{\GL}{\operatorname{GL}}
\newcommand{\PSL}{\operatorname{PSL}}

\newcommand{\SO}{\operatorname{SO}}
\newcommand{\OO}{\operatorname{O}}
\newcommand{\CO}{\operatorname{CO}}
\newcommand{\PO}{\operatorname{PO}}

\newcommand{\Ein}{\operatorname{\bf Ein}^{n-1,1}}
\newcommand{\tEin}{\widetilde{\operatorname{\bf Ein}}^{n-1,1}}
\newcommand{\Min}{\operatorname{\bf Min}^{n-1,1}}
\newcommand{\Euc}{\operatorname{\bf Euc}}

\newcommand{\Conf}{\operatorname{Conf}}

\newcommand{\BZ}{{\bf Z}}
\newcommand{\BN}{{\bf N}} 
\newcommand{\BR}{{\bf R}}
\newcommand{\BC}{{\bf C}}

\newcommand{\RP}{{\bf RP}}

\begin{document}

\title[Conformally flat Lorentzian manifolds]{Classification of closed conformally flat Lorentzian manifolds with unipotent holonomy}
\author{Rachel Lee and Karin Melnick}
\thanks{Melnick partially supported by NSF Award DMS-2109347; Lee partially supported by NSF Awards DMS-2109347 and DMS-2203493}
\maketitle

\linespread{1.1}
\section{Introduction}

A \emph{conformal structure} on a manifold $M$ is an equivalence class of semi-Riemannian metrics, where two metrics are equivalent if they are related by multiplication with a positive, smooth, real-valued function. A manifold that is locally conformally equivalent to a flat affine space is called \emph{conformally flat}. Such manifolds can alternatively be characterized as those admitting a $(G,{\bf X})$-structure where ${\bf X}$ is the suitable conformally flat, $G$-homogeneous model space. The locally homogeneous structure on $M$ gives rise to a developing pair $(\delta, \rho)$ where $\delta\colon \widetilde{M}\to {\bf X}$ is a local diffeomorphism from the universal cover $\widetilde{M}$ of $M$ to ${\bf X}$, and $\rho \colon \pi_1(M)\to G$ is the holonomy representation, such that $\delta$ is $\rho$-equivariant. Under group-theoretic assumptions on the holonomy image, classification results for $(G,{\bf X})$-manifolds can be obtained;  for example, W. Goldman proved:

\begin{theorem}[{\cite[Thm A]{goldman.nilp.holo}}]
Let $M$ be a closed, conformally flat, Riemannian manifold of dimension $n \geq 3$, and assume that the image of the holonomy representation is virtually nilpotent---that is, there exists a nilpotent subgroup of finite index. Then $M$ is finitely covered by the sphere ${\bf S}^n$, a flat torus ${\bf T}^n$, or a Hopf manifold, diffeomorphic to ${\bf S}^1\times {\bf S}^{n-1}$.
\end{theorem}

In each case, the developing map is proved to be a diffeomorphism onto an open subset of ${\bf X}$, in which case the structure is called \emph{Kleinian}, possibly all of ${\bf X}$, in which case it is \emph{complete}
(see Definition \ref{def.complete.kleinian} below). Then $M$ is geometrically isomorphic to the quotient of the developing image by the holonomy action.

This article concerns closed, conformally flat Lorentzian manifolds. We consider those with unipotent holonomy image. 
 The conformally flat model space in Lorentzian signature is $\Ein$.  Uniquely to this signature, it has infinite fundamental group, and the target of the developing map in general is the noncompact universal cover ${\bf X} = \tEin$.  These spaces are introduced in \S \ref{sec.prelims},
 see also the references \cite{frances.these, bcdgm.primer}.  Unipotence implies that the holonomy image stabilizes an isotropic flag in the standard representation of $G$, which corresponds to a chain of invariant subspaces in $\Ein$ 
\begin{equation}\label{eq1}
    \{p_0\} \subset \bar{\Delta} \subset L(p_0),
\end{equation}
consisting of a distinguished point, contained in a \emph{photon}, contained in the \emph{lightcone} of $p_0$.  The complement of a light cone in $\Ein$ is conformally equivalent to Minkowsi space $\mbox{\bf Min}^{n-1,1}$ and is called a \emph{Minkowski patch} (see \S\ref{sec.geom.ein} below).  Noting that $G$ acts transitively on isotropic flags, we fix one.  Our proof is organized according to the intersection of the developing image with the components of this flag. 

In conformal Riemannian geometry, the model is the round sphere, and the corresponding decomposition comprises just a point $p_0 \in {\bf S}^n$.  In that case, with unipotent holonomy, if the developing image contains $p_0$, then $M \cong {\bf S}^n$; otherwise, it is a flat torus.  
The Lorentzian case is considerably more complex.  Two cases are Lorentzian analogues of the Riemannian classification.  Between these are two intermediate cases, each giving rise to new examples.


%
%
The subspace $\bar{\Delta} \subset L(p_0)$ from \eqref{eq1} lifts in $\tEin$ to a photon $\Delta \subset L(\tilde{p}_0)$ contained in a light cone, both unbounded, for $\tilde{p}_0$ any point in the preimage of $p_0$.  The complement of $L(\tilde{p}_0)$ is a countable union of Minkowski patches.  We denote by $\widetilde{\mbox{O}}(2,n)$ the connected conformal group of $\tEin$; it is an infinite covering group of $\mbox{O}^0(2,n)$.  Its central element generating $\pi_1(\Ein) \cong \BZ$ will be denoted $\alpha$ (see \S \ref{sec.geom.tein} for details).  We denote by $\bar{\delta}$ the composition of $\delta$ with the covering $\pi_{\mathrm{Ein}} : \tEin \rightarrow \Ein$.

\begin{theorem}\label{mainthm} Let $M$ be a closed, conformally flat, Lorentzian manifold of dimension $n \geq 3$, with unipotent holonomy.  Then, up to composition of $\delta$ with a conformal equivalence of $\tEin$, one of the following holds:
\begin{enumerate}

    \item $p_0 \in \mbox{im } \bar{\delta}$: Then $M$ is a complete $(\widetilde{\mbox{O}}(n,2), \tEin)$-manifold.  The holonomy image is generated by an element of the form $\alpha^i g$ with $i\neq 0$ and $g$ projecting to a unipotent element of $\mbox{PO}^0(n,2)$. Topologically $M \cong {\bf S}^{n-1} \times {\bf S}^1$, up to a finite covering.
    
    \item $p_0 \notin \mbox{im } \bar{\delta}$ but $\mbox{im } \bar{\delta} \cap \bar{\Delta} \neq \emptyset$: Then $\delta$ is a diffeomorphism to a bounded open subset comprising the union of two Minkowski charts and the interstitial component of $L(\tilde{p}_0) \backslash \pi_{\mathrm{Ein}}^{-1}(p_0)$ lying in their common closure. The holonomy image is generated by one element 
which descends to act nontrivially on $\bar{\Delta}$.
    Topologically, $M \cong {\bf S}^{n-1}\times {\bf S}^1$ up to a finite covering.  
        
    \item $\mbox{im } \bar{\delta} \cap \bar{\Delta} =\emptyset$ but $\mbox{im } \bar{\delta} \cap L(p_0) \neq \emptyset$: Then $n = 2k+2$ for $k \in \BN$, and $\delta$ is a diffeomorphism to $\tEin \backslash \Delta$.  The holonomy image is a nilpotent extension by $\BZ$ of a discrete Heisenberg group of rank $2k+1$.  Topologically, $M$ is a nilmanifold, more precisely, a Heisenberg-nilmanifold bundle over ${\bf S}^1$ with unipotent monodromy.  
    
    \item $\mbox{im } \bar{\delta} \cap L(p_0) = \emptyset$:  Then $M$ is a complete $(O(n-1,1) \ltimes \BR^n, \mbox{\bf Min}^{n-1,1})$-manifold.  Topologically, $M$ is a nilmanifold of nilpotence degree at most 3.
      \end{enumerate}

In every case, $M$ is Kleinian.
\end{theorem}

Slightly more detailed statements and the proofs for cases (1), (2), (3), and (4) are in sections \ref{sec.case1}, \ref{sec.case2}, \ref{sec.case3}, and \ref{sec.case4.essential}, respectively.

One consequence of the above theorem is the following topological classification of closed conformally flat Lorentzian manifolds of dimension $n \geq 3$ with unipotent holonomy.
\begin{corollary}\label{maintheorem}
Let $M$ be a closed, conformally flat, Lorenztian manifold with unipotent holonomy, of dimension $n \geq 3$. Then $M$ is finitely covered by ${\bf S}^{n-1} \times {\bf S}^1$ or $M$ is a nilmanifold of degree at most 3.
\end{corollary}


One motivation for our classification is the goal of classifying closed Lorentzian manifolds admitting an \emph{essential} conformal flow---that is, a flow that does not preserve any metric in the conformal class.  In the Riemannian case, these can only be ${\bf S}^n$ by a celebrated theorem of Obata \cite{obata.lich} and Ferrand \cite{lf.lich}.  By the Lorentzian Lichnerowicz Conjecture, which has been proved for 3-dimensional, real-analytic Lorentzian manifolds \cite{fm.lich3d}, all essential examples should be conformally flat.  These can, however, be of infinitely-many different topological types (see \cite{frances.pq.counterexs}).
The conformally flat, closed Lorentzian manifolds admitting an essential conformal flow and having unipotent holonomy correspond to cases (1) and (2) in our classification theorem above.  In section \ref{sec.case4.essential} we obtain the following topological classification.

\begin{theorem}
\label{thm.essential}
Let $M$ be a closed, conformally flat Lorentzian manifold with unipotent holonomy of dimension $n \geq 3$.  Then the conformal structure on $M$ is essential if and only if $M$ is finitely covered by ${\bf S}^{n-1} \times {\bf S}^1$.
\end{theorem}

\subsection{Acknowledgments}

This paper generalizes to arbitrary dimension the dissertation result of the first author, in which she obtained the classification in dimension three.
We thank Bill Goldman for his role in proposing the problem, the Lorentzian analogue to Theorem 1.1.    As co-advisor of the first author's dissertation, his suggestion of using the techniques of \cite{fgh.nilp} for nilpotent holonomy lead to Proposition 2.16, which plays a key role in several places below.   We thank him for this and for further valuable feedback.

\section{Preliminary definitions and results}
\label{sec.prelims}

While conformally flat Riemannian manifolds are locally modeled on the round sphere with its group of M\"obius transformations, the \emph{Einstein universe} (also called the \emph{Lorentzian M\"obius space}) is the local model for conformally flat Lorentzian manifolds, and comes with a rank-two simple group of conformal transformations.  We put conformally flat Lorentzian geometry in the context of $(G,{\bf X})$-structures in the next subsection.  Then we present the geometry of the Einstein universe and its universal cover.  Finally, we focus on the action of the maximal connected unipotent subgroup of conformal transformations of the Einstein universe.

\subsection{The $(G,{\bf X})$-structure of a conformally flat Lorentzian manifold}
\label{sec.GX.structure}

Let $n \geq 3$.
Let $\mathcal{N} \subset \BR^{n+2}$ be the null cone of a nondegenerate quadratic form of index $2$.  Let ${\bf X} \subset {\bf RP}^{n+1}$ be the image of $\mathcal{N}$ under projectivization, a quadric hypersurface.  The restriction of the quadratic form to $T \mathcal{N}$ is a degenerate symmetric form which descends to a Lorentzian metric on ${\bf X}$, well-defined up to conformal equivalence; the resulting conformal Lorentzian manifold is $\Ein$. The group $\OO(2,n)$ of linear isometries of the quadratic form descends to a group of conformal transformations of $\Ein$ which is easily seen to be transitive.  The quotient $\PO(2,n)$ will be $G$.

Let $(M^n,g)$ be a conformally flat Lorentzian manifold.  For each $p \in M$, there is an open neighborhood $U$ of $p$ and a conformal diffeomorphism of $(U, \left. g \right|_U)$ with an open subset of $\Min$.  Minkowski space conformally embeds in $\Ein$, which is shown in \ref{sec.geom.ein} below.  In fact, $\Ein$ is the conformal completion of $\Min$ in the following sense:

\begin{theorem}(see \cite{frances.these} Thm 2.13, \cite{sharpe} Thm 5.2)
\label{thm.liouville}
Let $U,V \subset \Min$ be connected open subsets.  Fix a conformal embedding $\iota: \Min \rightarrow \Ein$.  Let $f : U \rightarrow V$ be a conformal diffeomorphism.  Then there is a unique $F \in G = \PO(2,n)$ such that $f$ is conjugate by $\iota$ to $\left. F \right|_{\iota(U)}$. 
\end{theorem}

This is the Lorentzian version of the Liouville Theorem.  
A consequence is the following Development Theorem.

\begin{theorem}(compare \cite[Thm 1.1]{goldman.nilp.holo})
Let $(M^n,g)$ be a conformally flat Lorentzian manifold with universal cover $\pi_M : \widetilde{M} \rightarrow M$.  Then there exists a pair $(\bar{\delta}, h)$ with
$\bar{\delta} : \widetilde{M} \rightarrow \Ein$ a conformal immersion and $h : \pi_1 (M) \rightarrow \PO(2,n)$ a homomorphism such that the diagram

\begin{equation*}
    \begin{tikzcd}
    \widetilde{M} \arrow[r, "\bar{\delta}"] \arrow[d, swap, "\gamma"]
    & \Ein \arrow[d, "h(\gamma)"] \\
    \widetilde{M} \arrow[r, "\bar{\delta}"]
    & \Ein
\end{tikzcd}
\end{equation*}
commutes for all $\gamma \in \pi_1 (M)$.  Moreover, if $(\bar{\delta}', h')$ is another such pair, then there exists $g \in \PO(2,n)$ such that $\bar{\delta} = g \circ \bar{\delta}'$ and $h'(\gamma) = g h (\gamma) g^{-1}$ for all $\gamma \in \pi_1(M)$.
\end{theorem}

More generally, we can take the existence of such a developing pair for any $(G,{\bf X})$ to be the definition of a \emph{$(G,{\bf X})$-structure} on $M$.  The following are standard terms.

\begin{definition}
\label{def.complete.kleinian}
\label{def_complete} Let $(\delta, h)$ be a developing pair for a $(G,{\bf X})$-structure on $M$.  Let $\Gamma < G$ be the image of $h$.  The $(G,{\bf X})$-structure is
\begin{enumerate}
    \item  \emph{complete} if $M \cong {\bf X} / \Gamma$
    \item \emph{Kleinian} if $M \cong \Omega/ \Gamma$ for $\Omega \subset {\bf X}$ an open subset.
\end{enumerate}
\end{definition}
In the case that ${\bf X}$ is not simply connected, the developing map always lifts to the universal cover $\widetilde{\bf X}$, and the holonomy lifts to $\widetilde{G}$, the group of lifts of $G$ to $\widetilde{\bf X}$.  One often prefers to speak of being complete or Kleinian with respect to the $(\widetilde{G}, \widetilde{\bf X})$-structure, but one may use both notions.

The following lemma will be applied to the developing map 
in order to conclude completeness in the sequel.

\begin{lemma}[see {\cite[Lem 3.4]{ch.daly.thesis}}]
\label{daly4}
Let $F: U \to X$ be a local diffeomorphism.  Let $W\subseteq U$ be an open subset on which $F$ restricts to a diffeomorphism onto $X$.  Assuming $U$ is connected, then $W = U$ and $F$ is a diffeomorphism.
\end{lemma}

We introduce here a few techniques for studying developing pairs for $(G,{\bf X})$-structures, which will be refined for our particular setting in subsequent sections.
The general idea is that holonomy-invariant objects on ${\bf X}$ correspond to well-defined objects on $M$.  Assuming $M$ is closed, these objects will provide leverage to establish completeness.  A first example is the following proposition, of which the short and easy proof is left to the reader.

\begin{proposition}
\label{prop.closed.invt}
Let $(\delta, h)$ be a developing pair for a $(G,{\bf X})$-structure on $M$, and let $\Gamma$ be the image of $h$.  Let $V \subset {\bf X}$ be closed and $\Gamma$-invariant.  Then $\pi_M(\delta^{-1}(V)) \subseteq M$ is closed.
\end{proposition}

Because $\delta$ is a local diffeomorphism, vector fields on ${\bf X}$ have well-defined pull-backs to $\widetilde{M}$.  In fact, the same is true for vector fields on submanifolds $V \subset {\bf X}$.  For $Y \in \mathcal{X}(V)$, the pull-back to $\delta^{-1}(V)$ will be denoted $\delta^*Y$.  

\begin{proposition}
\label{prop.vf.pullback}
Let $(\delta, h)$ be a developing pair for a $(G,{\bf X})$-structure on $M$, and let $\Gamma$ be the image of $h$.  Suppose that a regular submanifold $V \subset {\bf X}$ and a complete vector field $Y \in \mathcal{X}(V)$ are $\Gamma$-invariant.  Let $\hat{V}$ be a connected component of $\delta^{-1}(V)$.  If $\hat{V}$ is closed, then $\delta^*Y$ is complete, and
the image $\delta(\hat{V}) \subset V$ is invariant by the flow along $Y$.  
\end{proposition}

\begin{proof}
Let $\hat{\Gamma}$ be the group of deck transformations of $\widetilde{M}$.  By $\Gamma$-invariance of $Y$, the pullback $\delta^*Y$ is $\hat{\Gamma}$-invariant on $\delta^{-1}(V) = \hat{\Gamma}.\hat{V}$.  The latter is a union of connected components, each of which is closed in $\widetilde{M}$.  Therefore the image $\pi_M(\hat{V})$ is closed in $M$.  The vector field $\delta^*Y$ pushes forward to this image and the push-forward is complete.  Therefore $\delta^*Y$ is complete on $\hat{V}$.  By design, $\left. \delta \right|_{\hat{V}}$ intertwines the two flows, so $\delta(\hat{V})$ is invariant by the flow along $Y$.
\end{proof} 

\subsection{The geometry of the Einstein space and its universal cover}
\label{sec.geom.ein.tein}

This subsection details some of the analytic and synthetic geometry of $\Ein$ and $\tEin$.  Identities for causally defined sets are established, which will be used in the construction of examples in Section \ref{sec.case2} below.

\subsubsection{Geometry of $\Ein$}
\label{sec.geom.ein}

Recall that the construction of $\Ein$ begins with a nondegenerate, index-2 quadratic form on $\BR^{n+2}$.  It is convenient to fix the following one
$$ q_{n,2}(x) = 2 x_0 x_{n+1} + 2 x_1 x_n + \sum_{i=2}^{n-1} x_i^2$$
and to define for $x \in \BR^n$
$$ q_{n-1,1}(x) = 2 x_1 x_n + \sum_{i=2}^{n-1} x_i^2$$
which is of index 1.

Consider the following immersion of $\Min \rightarrow \BR^{n+2}$
$$ (x_1, \ldots, x_n) \mapsto (-q_{n-1,1}(x)/2, x_1, \ldots, x_n, 1)$$
This is a semi-Riemannian immersion of $\Min$ to $\BR^{n,2} = (\BR^{n+2}, q_{n,2})$.  The image is in the null cone $\mathcal{N}$ and is transverse to the fibers of the projectivization map.  Thus the composition
$$ \iota : (x_1, \ldots, x_n) \mapsto [-q_{n-1,1}(x)/2: x_1: \cdots : x_n : 1]$$
defines a conformal embedding of $\Min$ in $\Ein$, called a \emph{Minkowski chart}.  The image of such an embedding will also be called a \emph{Minkowski patch} below.  

The complement of the above Minkowski patch is the intersection of $\Ein$ with the subvariety of ${\bf RP}^{n+1}$ defined by $x_{n+1}=0$ in homogeneous coordinates.  According to $q_{n,2}$, this latter subvariety is the projectivization of $e_0^\perp$.  The intersection $e_0^\perp \cap \mathcal{N}$ is the union of the totally istropic planes containing $e_0$.  The projectivization of a totally isotropic plane in $\BR^{n,2}$ is called a \emph{photon} in $\Ein$.  The projectivization of $e_0^\perp \cap \mathcal{N}$ comprises all the photons of $\Ein$ passing through $[e_0]$.  It is a singular hypersurface called the \emph{lightcone} of $p_0 = [e_0]$, denoted $L(p_0)$.  This is thus the complement of our Minkowski patch, which evidently is determined by $p_0$ and will accordingly be denoted $\mbox{Min}(p_0) \subset \Ein$.

Note that $\OO(n,2)$ acts transitively on isotropic flags comprising an isotropic line, a totally isotropic plane, and a degenerate hypersurface, as above, and so $\PO(n,2)$ acts transitively on configurations $p_0 \subset \bar{\Delta} \subset L(p_0)$, where $\bar{\Delta}$ is a photon through $p_0$.

The photon $\bar{\Delta}$ can be identified with ${\bf RP}^1$, in a geometric sense.  The stabilizer in $\PO(n,2)$ of a totally isotropic plane is a subgroup isomorphic to $\PSL(2,\BR)$.  Thus $\bar{\Delta}$ inherits a $1$-dimensional real-projective structure from the geometry of $\Ein$, isomorphic to that of ${\bf RP}^1$.

Topologically, $\Ein$ is homeomorphic to ${\bf S}^{n-1} \times {\bf S}^1/\langle \sigma \rangle$, where $\sigma$ is the antipodal map on both factors.  The fundamental group of $\Ein$ is isomorphic to $\BZ$.  Moreover, the metric corresponding to $g_{{\bf S}^{n-1}} \oplus - g_{{\bf S}^1}$, where $g_{{\bf S}^k}$ is the constant-curvature metric on ${\bf S}^k$, belongs to the conformal class of $\Ein$.  For these facts we refer to \cite{frances.these}, \cite[Sec 4]{bcdgm.primer}.

Another way to see the topology of $\Ein$ is via the following useful projection.  Let $\bar{\Delta}$ be a photon, corresponding to the projectivization of the totally isotropic subspace $\mbox{span}\{ u,v \} \subset \BR^{n,2}$.  Let
\begin{eqnarray*}
 \rho_{\bar{\Delta}} & : & \Ein \backslash \bar{\Delta}  \rightarrow  \bar{\Delta}  \\
 & & \left[ x \right]  \mapsto  \left[ \langle x, v \rangle u - \langle x, u \rangle v \right]
 \end{eqnarray*}
where the inner product is the one determined by $q_{n,2}$.  It is easily checked that $\rho_{\bar{\Delta}}([x])$ is independent of the choice of basis $\{ u,v \}$ or the choice of $x$ representing $[x]$.  
It has a well-defined value in $\bar{\Delta}$, as $\mbox{span} \{ u,v \}$ is the inverse image of $\bar{\Delta}$ in $\mathcal{N}$ and is a maximal isotropic subspace.   For any $p \in \bar{\Delta}$, the fiber $\rho_{\bar{\Delta}}^{-1}(p)$ is $L(p) \backslash \bar{\Delta}$.  For $p\neq q$ both in $\bar{\Delta}$, the intersection $L(p) \cap L(q)$ is precisely $\bar{\Delta}$.  It follows that $\rho_{\bar{\Delta}}$ is a submersion, the fibers of which form a foliation by hypersurfaces, exhibiting $\Ein \backslash \bar{\Delta}$ as diffeomorphic to $\BR^{n-1} \times {\bf S}^1$.

\subsubsection{Geometry of $\tEin$}
\label{sec.geom.tein}

The universal covering $\tEin$ is homeomorphic to ${\bf S}^{n-1} \times\BR$ such that $\sigma$ lifts to 
\begin{equation*}
\begin{matrix}
    \tilde{\sigma}:  & {\bf S}^{n-1} \times\BR &\to & {\bf S}^{n-1} \times\BR \\
       & (x,t) &\mapsto &(-x, t+\pi)
\end{matrix}
\end{equation*}
The generator $\alpha$ of $\pi_1(\Ein)$ is represented by the deck transformation corresponding to $\tilde{\sigma}$ under this identification.  Each photon in this model is the graph of a unit-speed curve in ${\bf S}^{n-1}$.

Given a photon $\Delta \subset \tEin$, with $\bar{\Delta} = \pi_{\mathrm{Ein}} (\Delta)$, the map $\rho_{\bar{\Delta}}$ from the previous section lifts to $\tEin$:  the fibration of $\Ein\backslash \bar{\Delta}$ lifts to a foliation of $\tEin \backslash \Delta$ by closed hypersurfaces; then any lift of $\rho_{\bar{\Delta}}$ corresponds to the quotient map to the leaf space of this foliation, which is diffeomorphic to $\BR$.  We will define a specific lift which will in fact be a map $\rho_\Delta : \tEin \backslash \Delta \rightarrow \Delta$ in \S \ref{sec.tau} below.

Note that the geometric isomorphism $\bar{\Delta} \cong {\bf RP}^1$ lifts to a $(\mbox{PSL}(2,\BR),{\bf RP}^1)$-structure on $\Delta$, in which it is isomorphic to $\widetilde{\bf RP}^1$, with its transitive $\widetilde{SL}(2,\BR)$-action.

\begin{figure}
	\centering
    \includegraphics[]{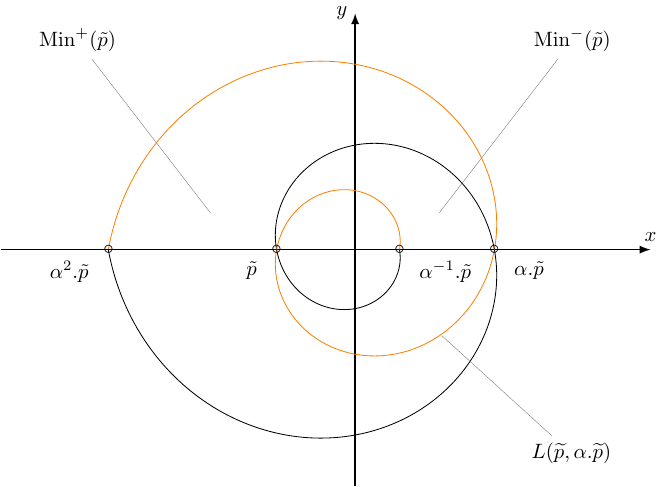}
	\caption{The above figure shows a cross-section of $\tEin \cong \BR^n\setminus\{0\}$, with a lightcone $L(\tilde{p})$ indicated by solid lines, and a photon in it indicated in orange.}
	\label{fig:Ein}
\end{figure}

In a useful refinement of the model for $\tEin$, we identify ${\bf S}^{n-1} \times \BR$ with $\BR^n\setminus\{0\}$ in the usual way, by
\begin{equation*}
    \begin{matrix}
        (x,\theta) &\mapsto & e^\theta x.
    \end{matrix}
\end{equation*}
Under this identification, $\alpha$ becomes
\begin{equation*}\label{alphaR}
\begin{matrix}
& z &\mapsto & -e^\pi z.    
\end{matrix}
\end{equation*}
Furthermore, each photon becomes a logarithmic spiral contained in a $2$-dimensional linear subspace of $\BR^n$.  The lightcone at a point $\tilde{p}$ is the revolution of any photon through $\tilde{p}$ around the axis connecting $\tilde{p}$ and $\alpha(\tilde{p})$. The complement of $L(\tilde{p})$ is a disjoint union of Minkowski patches, comprising the lifts of $\mbox{Min}(p)$, where $p = \pi_{\mathrm{Ein}} (\tilde{p})$.  The connected components are permuted by $\alpha$. 
Figure \ref{fig:Ein} shows a 2-dimensional cross-section of a light cone and a photon in it.  The Minkowski patches are the regions bounded between successive loops.
Two distinguished Minkowski patches adjacent to $\tilde{p}$
are labeled $\mbox{Min}^+(\tilde{p})$ and $\mbox{Min}^-(\tilde{p})$; these have a causal interpretation, which is given in the next section. 
A reference for this material is \cite[Sec 4.3]{bcdgm.primer}.
 
\subsubsection{Causal Geometry of $\tEin$}
\label{sec.causal.ein}

A reference for some of the basic material on causality is \cite[Chs 1--2]{penrose.diff.top.rel}.  
Recall that a \emph{causal} tangent vector is timelike or null, and a \emph{causal curve} is one with causal velocity.
%
A \emph{time-orientation} of a Lorentzian manifold is given by a timelike vector field, and a simply connected Lorentzian manifold is always time-orientable.  
A causal tangent vector is
\emph{future-pointing}, respectively \emph{past-pointing}, if its inner product with the time-orienting vector field is negative, respectively positive, and similarly for a causal curve. 
A piecewise smooth curve is called \emph{timelike}, \emph{lightlike}, or \emph{causal}, respectively, if the velocity vectors, including the velocity from above and below at break points, are of the corresponding type; moreover, at the break points both velocity vectors must have the same time orientation---that is, all velocity vectors along the curve are future-pointing or all are past-pointing.

\begin{definition}
For $M$ a time-orientable Lorentzian manifold, let $x,y\in M$.
\begin{enumerate}
    \item $x$ \emph{chronologically precedes} $y$ (often denoted $x\ll y$) if there is a future-directed timelike curve from $x$ to $y$.
    \item $x$ \emph{causally precedes} $y$ (often denoted $x \leq y$) if there is a future-directed causal curve from $x$ to $y$.
    \item $M$ is \emph{causal} if $x \leq y, y \leq x \iff x=y$.
\end{enumerate}
\end{definition}

If $x\leq y$ but $x\neq y$, we write $x < y$.

Under the conformal equivalence 
$$\tEin \cong ({\bf S}^{n-1} \times \BR, g_{\bf{S}^{n-1}} \oplus - d\theta^2)$$
 a time-orientation is given by the vector field $\partial_\theta$.
The map $\alpha$ preserves time-orientation;
thus $\Ein$ is time-orientable, although it is not orientable (see \cite[\S 4.2]{bcdgm.primer}).
It is also causal (in fact, it is \emph{globally hyperbolic}). 
For future reference, a natural choice of coordinate $\theta$ on $\Ein$ gives, in homogeneous coordinates
\begin{equation}
\label{eqn.dtheta}
\partial_\theta   =   (x_0 - x_{n+1}) \left( \partial_1 - \partial_n \right) + (x_n - x_1) \left(  \partial_0 -  \partial_{n+1} \right)
\end{equation}

\begin{definition} Let $S\subset M$. The \emph{causal future set} $J^+(S)$ and \emph{future set} $I^+(S)$ are defined as
\begin{equation*}
    J^+(S)=\{y\in M:\exists x\in S, \,  x\leq y\} \quad
    I^+(S)=\{y\in M:\exists x\in S, \,  x\ll y\}.
\end{equation*}
The \emph{causal past set} and the \emph{past set} are defined analogously. 
\end{definition}

The future and past sets of a point in $\tEin$ have the following identifications in ${\bf S}^{n-1} \times\BR \cong \tEin$.
\begin{lemma}
\label{lem_casonEin}
Let $p=(x_0,t_0) \in \tEin$. Then the causal past and future can be expressed as
\begin{equation*}
    J^-(p) = \{(x,t) \in {\bf S}^{n-1} \times \BR 
    : t_0-t \geq d(x,x_0)\}, \quad J^+(p) = \{(x,t) \in {\bf S}^{n-1} \times \BR 
    : t-t_0 \geq d(x,x_0)\} 
\end{equation*}
and the past $I^-(p)$ and the future $I^+(p)$ as
\begin{equation*}
    I^-(p) = \{(x,t) \in {\bf S}^{n-1} \times \BR  
    : t_0-t > d(x,x_0)\}, \quad I^+(p) = \{(x,t) \in {\bf S}^{n-1} \times \BR  
    : t-t_0 > d(x,x_0)\} 
\end{equation*}
where $d$ denotes the standard Riemannian distance on ${\bf S}^{n-1}$.
\end{lemma}

The proof is left to the reader.

\begin{remark}
An immediate consequence of the above lemma is that $x\in\tEin$ always causally precedes $\alpha x$.
\end{remark} 

%
%
The causal future, respectively past, sets are the closure of the future, respectively past, sets. 
%
%
Note that given a point $p$, its future set and its past set are not complementary; however, 
\begin{lemma}\label{Eincomplement}
For any point $p = (x_0, t_0)\in \tEin$,
\begin{equation*}
    \tEin\setminus I^+(p) = J^-(\alpha .p)
\end{equation*}
\begin{equation*}
    \tEin\setminus J^+(p) = I^-(\alpha. p)
\end{equation*}
\end{lemma}
\begin{proof}
We prove only the first statement; the second is proved \emph{mutatis mutandis}. From Lemma \ref{lem_casonEin},
\begin{equation*}
    J^-(\alpha p) = \{(x,t): (t_0+\pi)-t \geq d(x, -x_0)\}
\end{equation*}
As
$d(-x_0, x) + d(x, x_0) = \pi$ for any $x\in {\bf S}^{n-1}$, we have, again using Lemma \ref{lem_casonEin},
\begin{align*}
    \tEin\setminus I^+(p)
    &= \{(x,t): t-t_0 \leq d(x, x_0)\}\\
    &= \{(x,t): t-t_0 \leq \pi - d(x, -x_0)\}\\
    &= \{(x,t): t-(t_0 +\pi) \leq - d(x, -x_0)\}\\
    &= \{(x,t): (t_0 +\pi)-t \geq  d(x, -x_0)\}\\
    &= J^-(\alpha p),
\end{align*}
as desired.
\end{proof}

For $p \in \tEin$, the distinguished Minkowski patch $\mbox{Min}^+(p)$ is the Minkowski patch containing $p$ in its closure and in the future of $p$, while $\mbox{Min}^-(p)$ is the Minkowski patch containing $p$ in its closure, the points of which are not causally related to $p$.  We will take as definition the following, and verify a little later that these project to a Minkowski patch in $\Ein$ as previously defined.


\begin{definition}\label{Mincausal} For $p \in \tEin$,
\begin{enumerate}
    \item $\mbox{\bf Min}^+(p) = I^+(p)\cap I^-(\alpha^2 .p)$,
    \item $\mbox{\bf Min}^-(p) = I^+(\alpha^{-1}.p)\cap I^-(\alpha .p)$.
\end{enumerate}
\end{definition}

The following relation follows immediately from Definition \ref{Mincausal}.
\begin{proposition}
\label{cor.min.pm}
For any $p \in \tEin$,
$$\mbox{\bf Min}^+(p) = \mbox{\bf Min}^-(\alpha. p)$$
\end{proposition}


For two points $p$ and $q$ such that $q = \alpha^i .p$ for some $i\in\BZ_{>0}$, we define the following subsets, contained in $L = L(p) \cap L(q)$:
\begin{equation*}
\begin{matrix}
    L(p,q)=\{x\in L: p <  x < q\}, & 
    L[p,q]=\{x\in L: p \leq  x \leq q\}\\
  \end{matrix}
\end{equation*}

The boundaries of the Minkowski patches are the following subsets of $L(p)$.
\begin{corollary}\label{Minboundary} For any $p \in \tEin$,
$$\partial\mbox{\bf Min}^+(p)=L[p,\alpha^2. p] \qquad \partial\mbox{\bf Min}^-(p)=L[\alpha^{-1}.p, \alpha .p]$$
\end{corollary}

\begin{proof}
For $p=(\phi_0,\theta_0) \in {\bf S}^{n-1} \times \BR$, Definition \ref{Mincausal} and Lemma \ref{lem_casonEin} give 
\begin{align*}
    \mbox{\bf Min}^+(p) 
    &= \{(\phi,\theta) \in {\bf S}^{n-1} \times \BR :  d(\phi, \phi_0) < \theta - \theta_0 < 2 \pi - d(\phi, \phi_0) \}
\end{align*}
Using that $0\leq d(\phi, \phi_0)\leq \pi$ for all $\phi$,
the boundary set
\begin{equation*}
   \{ \theta-\theta_0=d(\phi,\phi_0)\} = L[p,\alpha. p]
    \end{equation*}
    while
   \begin{equation*}
    \{ (\theta_0+2\pi)-\theta=d(\phi_0,\phi)\} = L[\alpha .p, \alpha^2 .p]
\end{equation*}

This proves the equality for $\partial \mbox{\bf Min}^+(p)$; the second follows from this one and Corollary \ref{cor.min.pm}.
%
\end{proof}

For a fixed $p = (\phi_0, \theta_0)$ and given any $q = (\phi, \theta)$, there is $k \in \BZ$ such that
$$ d(\phi, \phi_0) \leq \theta - \theta_0 - 2k \pi \leq 2 \pi - d(\phi, \phi_0) \ \mbox{or} \  - d(\phi, \phi_0) \leq \theta - \theta_0 - 2k \pi \leq d(\phi, \phi_0)$$

Strict inequality in the first case corresponds to $q \in \mbox{\bf Min}^+(\alpha^{2k}.p)$, while strict inequality in the second case corresponds to $q \in \mbox{\bf Min}^-(\alpha^{2k}.p) = \mbox{\bf Min}^+(\alpha^{2k-1}.p)$.  Equality corresponds to $p \in L[\alpha^{2k}.p,\alpha^{2k+2}.p]$ or $p \in L[\alpha^{2k-1}.p,\alpha^{2k+1}.p]$, respectively.  Thus $\tEin$ is the disjoint union over $i \in \BZ$ of the Minkowski patches $\mbox{\bf Min}^+(\alpha^i.p)$ with $L(p)$.  The first union and $L(p)$ are each $\alpha$-invariant.  It follows that the projection $\pi_{\mathrm{Ein}}(\mbox{\bf Min}^+(p))$ is the complement of $L(\bar{p})$ in $\Ein$, which is the Minkowski patch $\mbox{\bf Min}(\bar{p})$.  The projection of $\mbox{\bf Min}^-(p)$ also equals $\mbox{\bf Min}(\bar{p})$.

%
%
%
\subsection{Unipotent dynamics on $\Ein$}
\label{sec.dyn.ein}

\subsubsection{The maximal unipotent subgroup}
\label{sec.unip.subgp}

The maximal unipotent subgroup of $\hat{G}=\OO(n,2)$ is the unipotent radical of the minimal parabolic subgroup.  (These groups are of course unique only up to conjugacy in $\hat{G}$.) The latter subgroup is the stabilizer of an isotropic flag of $\BR^{n,2}$; the unipotent is the subgroup of the stabilizer of an isotropic flag which moreover acts trivially on each subquotient of the flag.  For the quadratic form $q_{n,2}$ and the isotropic flag
$$ \mathcal{F} = \BR e_0  \subset \mbox{span} \{ e_0, e_1 \} \subset e_0^\perp \subset \BR^{n,2}$$
the maximal unipotent subgroup is upper-triangular.  It will be denoted $\mathcal{U}$.  This group is simply connected, and is contained in $\hat{G}^0$.
As $-\mbox{Id}_{n+2}$ is not in $\mathcal{U}$, it projects isomorphically to its image in $G= \PO(n,2)$, which we will also denote by $\mathcal{U}$.

Recall that $\widetilde{\OO}(n,2)$ denotes the conformal group of $\tEin$, for which we may also write $\widetilde{G}$.  (Note that $\pi_1(\OO^0(n,2)) \cong \BZ_2 \times \BZ$ for $n \geq 3$, so $\widetilde{G}^0$
is a two-fold quotient of the universal cover of $\OO^0(n,2)$.)   Denote by $q : \widetilde{\OO}(n,2) \rightarrow \PO(n,2)$ the quotient.  As $\pi_1(\Ein) \cong \BZ$ this is a $\BZ$-covering on the identity components.  We will denote by $\widetilde{\mathcal{U}}$ the full lift $q^{-1}(\mathcal{U})$ and by $\widetilde{\mathcal{U}}^0$ the identity component; the latter group is also isomorphic to $\mathcal{U}$, as $\mathcal{U}$ is simply connected.

The isotropic flag $\mathcal{F}$ stabilized by $\mathcal{U}$ corresponds to the chain of subspaces
$$ p_0 = [e_0] \subset \bar{\Delta} \subset L(p_0) \subset \Ein$$
where $\bar{\Delta}$ is the photon obtained from the projectivization of $\mbox{span} \{ e_0,e_1 \}$.  
The $\mathcal{U}$-action on $\bar{\Delta} \cong {\bf RP}^1$ is the projective parabolic flow fixing $p_0$.  The $\widetilde{\mathcal{U}}^0$-action on $\Delta \subset \tEin$ corresponds under the geometric isomorphism $\Delta \cong \widetilde{\bf RP}^1$ to the lift of this one-parameter subgroup to $\widetilde{\SL(2,\BR)}$.

Because it stabilizes the complementary light cone, $\mathcal{U}$ acts conformally on the Minkowski patch $\mbox{\bf Min}(p_0)$.  By the Liouville Theorem \ref{thm.liouville}, this action is faithful.  The conformal group of $\Min$ is the similarity group $\CO(n-1,1) \ltimes \BR^n$.  For $g \in \mathcal{U}$, we denote by $L_g$ the image under the homomorphism $\mathcal{U} \rightarrow \CO(n-1,1)$.  The image $L(\mathcal{U})$ is a unipotent subgroup of $\SO(n-1,1)$, stabilizing the isotropic flag $\BR e_1 \subset e_1^\perp \subset \BR^{n-1,1}$ for the quadratic form $q_{n-1,1}$.  We denote by $u_g$ the $\BR^n$-component of $g \in \mathcal{U}$, so that $g$ acts on $\Min(p_0) \cong \Min$ by the affine transformation $v \mapsto L_g(v) + u_g$.  Note that $u : \mathcal{U} \rightarrow \BR^n$ is a 1-cocycle for the representation of $\mathcal{U}$ on $\BR^n$ via $L$.

\subsubsection{The $\tau$-flow}
\label{sec.tau}

The maximal unipotent subgroup $\mathcal{U}$ has one-dimensional center that corresponds to a translation by the isotropic vector $e_1$ on $\Min$.
We will refer to the $\BR$-action of $Z(\mathcal{U})$ on $\Ein$ as the \emph{$\tau$-flow}.  In homogeneous coordinates on $\Ein$, it is
\begin{equation}
\label{eqn.tau.flow}
 \tau^s. [ x_0 : \cdots : x_{n+1}] = [x_0 + s x_n : x_1 - s x_{n+1} : x_2 : \cdots : x_{n+1}]
 \end{equation}
The lift to $\widetilde{\mathcal{U}}^0$ acting on $\tEin$ will be denoted by $\tau$ as well.  We will denote the corresponding vector fields on $\Ein$ and $\tEin$ by $Y_\tau$. In homogeneous coordinates on $\Ein$ it is $Y_\tau = x_n \partial_0 - x_{n+1} \partial_1$.  Using the time orientation $\partial_\theta$ from (\ref{eqn.dtheta}), the inner product $\langle Y_\tau, \partial_\theta \rangle \leq 0$, vanishing exactly on $\bar{\Delta}$.  Thus $Y_\tau$ is future-pointing on $\Ein \backslash \bar{\Delta}$ and $\tEin \backslash \Delta$.

Each point of $\Ein \backslash \bar{\Delta}$ tends under $\tau^s$ to a point of the fixed set $\bar{\Delta}$ as $s \rightarrow \pm \infty$.  The limit point is given by the projection $\rho_{\bar{\Delta}}$ from \S \ref{sec.geom.ein}.  
As a consequence, we have the following description, for $q \in \bar{\Delta}$:
\begin{equation*}
    \mbox{\bf Min}(q) = \{p\in\Ein\setminus \bar{\Delta}: \lim_{t\to\pm\infty}\tau^t.p = \rho_{\bar{\Delta}}(p) \neq q\}.
\end{equation*}

Similarly, the fixed set in $\tEin$ of $\{ \tau^s \}$ is $\Delta$ and points of $\tEin \backslash \Delta$ tend under $\tau^s$ to a point of $\Delta$ as $s \rightarrow \pm \infty$.
 Each Minkowski patch in $\tEin\setminus \Delta$ can be written as a set of points converging to a particular segment of $\Delta$ under the $\tau$-flow. 
For any $p, q \in \Delta$, denote $\Delta[p,q] = L[p,q] \cap \Delta$, and similarly for $\Delta(p,q)$.
\begin{proposition}\label{prop.minlimit}
For any $p\in\Delta,$ 
\begin{equation*}
\mbox{\bf Min}^\pm(p) = \{{x}\in\tEin\setminus\Delta: \lim_{t\to\mp\infty}\tau^t. {x} \in \Delta(p,\alpha. p)\}.
\end{equation*}
\end{proposition}

\begin{proof}
From Corollary \ref{Minboundary}, we know that
\begin{equation*}
    \Delta\cap \overline{\mbox{\bf Min}^+({p})}
    =\Delta[{p},\alpha^2 .p].
\end{equation*}

For ${x}\in\mbox{\bf Min}^+({p})$, let
$\bar{x}=\pi_{\mathrm{Ein}}({x})$ and $\bar{p}=\pi_{\mathrm{Ein}}(p)$. 
Let 
$$\bar{q} = \lim_{t \rightarrow \infty} \tau^t .\bar{x}  =  \lim_{t \rightarrow - \infty} \tau^t .\bar{x} \in\bar{\Delta}\setminus\{\bar{p}\}$$
and let ${q}$ be the $\pi_{\mathrm{Ein}}$-preimage of $\bar{q}$ in $\Delta({p}, \alpha {p})$.
Since ${x}\notin L({p}) = L(\alpha^i .p)$ for all $i$, we have
\begin{equation*}
    \mbox{\bf Min}^+({p})\subseteq \{{x}\in\tEin\setminus \Delta : \lim_{t\to-\infty}\tau^t. {x} \in \Delta({p}, \alpha .p)\cup \Delta(\alpha .p, \alpha^2 .p)\}.
\end{equation*}
Thus
$\displaystyle\lim_{t\to\pm\infty}\tau^t.{x}\in \{ {q}, \alpha.{q} \}$.  
Note that the path $t \mapsto \tau^t.{x}$ for $t \in (- \infty, \infty)$ is future-pointing, as $Y_\tau$ is future-pointing. If the two limits $\displaystyle\lim_{t\to\pm\infty}\tau^t.{x}$ were equal, the result would be a closed causal loop, violating causality of $\tEin$ (see \S \ref{sec.causal.ein}).  Moreover, the forward limit must be in the future of the backward limit, from which we conclude
\begin{equation*}
    \lim_{t\to-\infty}\tau^t.{x} = {q}\in \Delta({p}, \alpha. {p}) \qquad
    \lim_{t\to\infty}\tau^t.{x} = {q}\in\alpha(\Delta({p}, \alpha.{p}))
\end{equation*}

On the other hand, every point of $\bar{\Delta}\setminus \{ \bar{p} \}$ is the forward and backward limit of a point of $\mbox{\bf Min}(\bar{p})$ under $\{ \tau^t \}$.  As $\alpha$ commutes with $\tau^t$, every ${q} \in \Delta({p}, \alpha .{p})$ is $\lim_{t \rightarrow - \infty} \tau^t.{x}$ for some ${x} \in \mbox{\bf Min}^+(\alpha^i. {p})$ for some $i$; also every ${q} \in \Delta(\alpha. {p}, \alpha^2 .{p})$ is $\lim_{t \rightarrow  \infty} \tau^t.{x}$ for some ${x} \in \mbox{\bf Min}^+(\alpha^j .{p})$ for some $j$.  By the inclusions proved in the previous paragraph, necessarily $i=j=0$, and 
\begin{eqnarray*}
\mbox{\bf Min}^+({p}) & = & \{{x}\in\tEin\setminus\Delta: \lim_{t\to -\infty}\tau^t. {x} \in \Delta({p},\alpha .{p})\} \\
& = & \{{x}\in\tEin\setminus\Delta: \lim_{t\to \infty}\tau^t. {x} \in \Delta(\alpha. {p},\alpha^2 .{p})\} 
\end{eqnarray*}
The desired identity for $\mbox{\bf Min}^-({p})$ now follows from Corollary \ref{cor.min.pm}.
\end{proof}

We now define $\rho_\Delta : \tEin \backslash \Delta \rightarrow \Delta$.  Let $x \in \tEin$ with $\bar{x} = \pi_{\mathrm{Ein}}(x)$.

$$ \rho_\Delta(x) = \lim_{t \rightarrow \infty} \tau^t.x = \left\{  \begin{array}{cl}  \alpha^i .p & x \in L(\alpha^{i-1}.p,\alpha^i .p) \backslash \Delta \\
\pi_{\mathrm{Ein}}^{-1} (\rho_{\bar{\Delta}}(\bar{x})) \cap \Delta(\alpha^i.p,\alpha^{i+1}.p) & x \in \mbox{\bf Min}^-(\alpha^i.p)
\end{array}
\right.
$$

This map is a submersion onto $\Delta$ lifting $\rho_{\bar{\Delta}}$.

\subsubsection{An iterative technique for establishing completeness}

The following is a general technique for showing that the developing map is a diffeomorphism onto eligible subsets of the model space ${\bf X}$.  The basic version appears for affine manifolds with unipotent holonomy in the proof of completeness of \cite[Thm 6.8]{fgh.nilp} and our proof idea is derived from theirs.  

\begin{proposition}
\label{prop.fgh.technique}
Let $(\delta,h)$ be a developing pair for a $(G,{\bf X})$-structure on a closed manifold $M$, with holonomy image $\Gamma < G$.  Let $V \subset {\bf X}$ be a connected, $\Gamma$-invariant regular submanifold, with a $\Gamma$-invariant foliation $\mathcal{F}$.  Let $\hat{V}$ be a connected component of $\delta^{-1}(V)$.  Assume $\hat{V}$ is a closed set.

\begin{enumerate}
\item{
Suppose there are complete vector fields $Y_1, \ldots, Y_d \in \mathcal{X}(V)$, such that, for all $i$,
 \begin{itemize}
 \item $Y_i$ is locally projectable modulo $\mathcal{F}$
 \item $\gamma_* Y_i \equiv Y_i \ \mbox{mod } T \mathcal{F} \ \forall \gamma \in \Gamma$
 \end{itemize}
 If the image $\delta(\hat{V})$ is $\mathcal{F}$-saturated, then $\delta(\hat{V}) \subset V$ is invariant by the flow along every $Y_i$.}

\item{
Let 
$W \subseteq V$ be a connected, regular submanifold saturated by $\mathcal{F}$ and by the flows along $\{ Y_i \}$, such that
\begin{itemize}
\item at all $y \in W$, the projections mod $\mathcal{F}$ of $\{ Y_i \}$ form a frame of the local leaf space 
 \item $[Y_i, Y_j] \equiv 0 \ \mbox{mod } T \mathcal{F}$ on $W \ \forall i, j$ 
 \end{itemize}
  Denote by $\hat{\mathcal{F}}$ the pulled-back foliation of $\hat{V}$, 
  and let $\hat{W}$ be a connected component of $\delta^{-1}(W)$.  If  $\delta$ maps each leaf of $\hat{\mathcal{F}}$ in $\hat{W}$ diffeomorphically to its image,
 then $\left. \delta \right|_{\hat{W}}$
 is a covering map onto $W$.}
 \end{enumerate}
\end{proposition}

A vector field $Y$ is \emph{locally projectable} modulo a foliation $\mathcal{F}$ if every point has a foliated neighborhood $U$ with local leaf space $L$ such that for the projection $\rho_L : U \rightarrow L$, the push-foward $(\rho_L)_*Y$ is well-defined.

\begin{proof}
Begin with the assumptions of (1).  The first step is to build vector fields $\{ X_i \}$ on $\widetilde{M}$ corresponding to $\{ Y_i \}$.  
Denote by $\hat{\Gamma}$ the group of deck transformations of $\widetilde{M}$.
The saturation $\hat{\Gamma}.\hat{V}$ is a union of closed connected components, so its projection to $M$ is closed, as is the component $\pi_M(\hat{V})$.  For any $\bar{x} \in \pi_M(\hat{V})$, there is a neighborhood $\bar{U}$ and a diffeomorphism from $\bar{U}$ to a neighborhood $U \subset V$.  Let $\{ \bar{U}_j \}$ be a finite cover of $\pi_M(\hat{V})$ by such neighborhoods, and define vector fields $\bar{X}_1^j, \ldots, \bar{X}_d^j$ by pulling back $Y_1, \ldots, Y_d$ from $U \subset V$ to $\bar{U}_j$, for each $j$.  Let $\{ \psi_j \}$ be a partition of unity subordinate to $\{ \bar{U}_j \}$ and define $\bar{X}_i = \sum_j \psi_j \cdot \bar{X}_i^j$.  Let $X_i$ be the lift of $\bar{X}_i$ to $\hat{V} \subset \widetilde{M}$ for each $i$; it is complete.

By the assumption of $\Gamma$-invariance of $\mathcal{F}$, the pulled-back foliation $\hat{\mathcal{F}}$ descends to $\pi_M(\hat{V})$. For any $i$, by the $\Gamma$-invariance of $Y_i$ mod $\mathcal{F}$, 
on the overlap of $\bar{U}_j$ with some $\bar{U}_k$, the vector fields $\bar{X}_i^j \equiv \bar{X}_i^k$ modulo the foliation.  For the lifted vector fields $\{ X_i \}$, by construction, 
for any $x \in \hat{V}$ with neighborhood $\hat{U} \subset \hat{V}$ mapping diffeomorphically under $\delta$ to $U \subset V$, the push-forward $(\left. \delta \right|_{\hat{U}})_* X_i \equiv Y_i \ \mbox{mod } T \mathcal{F}$ for all $i$.  In particular, the $\{ X_i \}$ are projectable modulo $\hat{\mathcal{F}}$.

We proceed to prove (1).
Let $y = \delta(x)$ for $x \in \hat{V}$, and let $Y = Y_i$ for some $i \in \{ 1, \ldots, d \}$.  
Consider for arbitrary $t_0 > 0$ a path 
$\alpha(t) = \varphi^t_Y.y$, $ 0 \leq t \leq t_0$.
Let $X = X_i$ and let $\hat{\alpha}(t) = \varphi^t_X.x$ for $0 \leq t \leq t_0$, which is defined because $X$ is complete.
Let $\beta(t) = \delta \circ \hat{\alpha}(t)$.  

Now let $\tau = \{ t : \beta(t) \equiv \alpha(t) \}$ where equivalence means belonging to the same leaf of $\mathcal{F}$.  
By construction $0 \in \tau$.  Suppose that $t_k \rightarrow t$ with $t_k \in \tau$. By continuity of $\alpha$, $\beta$, and $\mathcal{F}$, the leaves
$$ F_{\alpha(t_k)} \rightarrow F_{\alpha(t)} \qquad \mbox{and} \qquad F_{\beta(t_k)} \rightarrow F_{\beta(t)} $$
so 
$$ F_{\alpha(t_k) } = F_{\beta(t_k)} \ \forall k \Rightarrow F_{\alpha(t)} = F_{\beta(t)}$$
Then $\tau$ is closed.  

For $t \in \tau$, let $L$ be a sufficiently small transversal so that the $\mathcal{F}$-holonomy $H_{\beta(t),\alpha(t)}$ is defined on $L$.  Denote by $\rho_L$ the projection from a small neighborhood $U$ of $\beta(t)$ to $L$.  Then $(\rho_L)_*Y$ is well-defined, by projectability.  Shrinking $L$ if necessary,  $(\rho_L)_*\delta_* X$ is also defined.  These projected vector fields are equal, because $\delta_* X \equiv Y$ modulo $\mathcal{F}$ wherever both are defined.  The projection $\rho_L \circ \beta$, where defined, is the integral curve through $\rho_L (\beta(t))$ of this vector field.  Now let $L' = H_{\beta(t),\alpha(t)}(L)$, a transversal to $\mathcal{F}$ through $\alpha(t)$.  Here $(\rho_{L'})_*Y$ is defined, and $\rho_{L'} \circ \alpha$, where defined, is its integral curve through $\rho_{L'}(\alpha(t))$.  Because $H_{\beta(t),\alpha(t)}$ pushes forward $(\rho_L)_*Y$ to $(\rho_{L'})_* Y$, it follows that $\rho_{L'} \circ \alpha = H_{\beta(t),\alpha(t)} \circ \rho_L \circ \beta$.  Then on a small interval of time around $t$ where these projections are both defined, $\alpha \equiv \beta$.  Thus $\tau$ is open.  We conclude $\tau = [0,t_0]$, so $\alpha(t_0) \equiv \beta(t_0) \in \delta(\hat{V})$.  By the hypothesis that $\delta(\hat{V})$ is $\mathcal{F}$-saturated, $\alpha(t_0) \in \delta(\hat{V})$, as desired.



Now let $W$ and $\hat{W}$ be as in (2).  Given $y \in W$, let $U \subset W$ be a foliated neighborhood with projection $L$ to the local leaf space.
The hypotheses imply that after possibly shrinking $U$, there is $\epsilon > 0$ such that the map
$$ B_{2 \epsilon}(\BR^d) \rightarrow L \qquad c = (c_1, \ldots, c_d) \mapsto \rho_L(\varphi^{c_d}_{Y_d} \circ \cdots \circ \varphi^{c_1}_{Y_1}.y)$$
is a diffeomorphism.  Let $L'$ be the image of $B = B_\epsilon(\BR^d)$.  For a point $y'$ of the leaf space, denote by $F_{y'} \subset W$ the corresponding leaf.  Let $U' = \cup_{y' \in L'} E_{y'}$ where $E_{y'} = U \cap F_{y'}$.  Given $c \in \BR^d$, denote by $\varphi^{\bf c}_{\bf Y}$ the composition $\varphi^{c_d}_{Y_d} \circ \cdots \circ \varphi^{c_1}_{Y_1}$.  

Given $c \in \BR^d$, let $\varphi^{\bf c}_{\bf X}$ denote the composition $\varphi^{c_d}_{X_d} \circ \cdots \circ \varphi^{c_1}_{X_1}$ on $\hat{V}$.  
Given $x \in \delta^{-1}(y) \cap \hat{W}$, reduce $\epsilon$ if necessary so that the map 
$$ c = (c_1, \ldots, c_d) \mapsto \rho_{\hat{L}} \circ \varphi^{\bf c}_{\bf X}(x)$$
is a diffeomorphism from $B$ to the local leaf space $\hat{L}$ of a foliated neighborhood of $x$.  We will assume that $\hat{L}$ is the full image of this map.
Then $\hat{L}$ maps bijectively via $\bar{d} = \rho_L \circ \delta \circ \rho_{\hat{L}}^{-1}$ onto $L'$.  Indeed, suppose $u = \rho_{\hat{L}} \circ \varphi^{\bf c}_{\bf X}.x$ and $u' = \rho_{\hat{L}} \circ \varphi^{\bf c'}_{\bf X}.x$ map under $\bar{d}$
to the same point in $L'$.  By equivariance of $\delta$ with respect to the flows and the foliations, that would mean $\rho_L (\varphi^{\bf c}_{\bf Y}.y) = \rho_L ( \varphi^{\bf c'}_{\bf Y}.y)$.  Then $\rho_L ( \varphi^{{\bf c } - {\bf c'}}_{\bf Y}.y) = \rho_L(y)$ because the $Y_i$s commute modulo $\mathcal{F}$.  Since $c - c' \in B_{2 \epsilon}(\BR^d)$, 
necessarily $ c =  c'$.  Thus injectivity is proved.  Surjectivity follows from equivariance of $\delta$, as well.  

For $u \in \hat{L}$, define $\hat{E}_u = (\left. \delta \right|_{\hat{F}_u})^{-1}(E_{\bar{d}(u)})$, where $\hat{F}_u$ is the $\hat{\mathcal{F}}$-leaf of $u$ in $\hat{V}$.  By assumption, these map diffeomorphically to their images under $\delta$ for all $u \in \hat{L}$. Thus
$ \hat{U}_x = \bigcup_{u \in \hat{L}} \hat{E}_u \subset \hat{W}$
maps bijectively, and hence diffeomorphically, to $U' \subset W$ under $\delta$.

Finally, suppose that for $x, x' \in \delta^{-1}(y)$, there is $u \in \hat{U}_x \cap \hat{U}_{x'}$.
Then $\hat{F}_{\rho_{\hat{L}}(u)} = \hat{F}_{\rho_{\hat{L}'}(u)}$, where $\hat{L}'$ is the local leaf space of $\hat{U}_{x'}$.
By construction, there exist unique $c, c' \in B$ such that $\varphi^{\bf c}_{\bf X}.x$ and $\varphi^{\bf c'}_{\bf X}.x'$ each belong to this leaf. Then $\varphi^{\bf - c'}_{\bf X}  \circ \varphi^{\bf c}_{\bf X}.x$ is in the same leaf as $x'$.  Then $\varphi^{\bf - c'}_{\bf Y}  \circ \varphi^{\bf c}_{\bf Y}.y$ is in the same leaf as $\delta(x') = y$.  As $c - c' \in B_{2 \epsilon}(\BR^d)$, the first point is in $U$, and the two points have the same projection under $\rho_L$.  By the assumption on $\epsilon$, it follows that $c=c'$.  But then $x$ and $x'$ are in the same leaf, which implies they are equal because leaves map diffeomorphically under $\delta$.  Therefore the open sets $\hat{U}_x$ are disjoint for distinct $x \in \delta^{-1}(y) \cap \hat{W}$.  

We conclude that $\delta$ is a covering from $\hat{W}$ onto its image in $W$.  The image is open.  If there were $y \in \partial (\delta(\hat{W}))$, then by the framing assumption of (2) it would be connected by a finite sequence of flows along $\{Y_i \}$ and segments in leaves of $\mathcal{F}$ to a point of $\delta(\hat{W})$.  But then by the flow-invariance of $\hat{W}$ and the assumption on $\delta$ along the foliations, $y$ would also be in the image $\delta(\hat{W})$.  We conclude $\delta$ is a covering map of $\hat{W}$ onto $W$, as desired.
\end{proof}


\section{The developing image contains $p_0$.}
\label{sec.case1}

Recall the $\mathcal{U}$-invariant flag 
$$p_0 \subset \bar{\Delta} \subset L(p_0) \subset \Ein$$

Our first case is when the image of $\bar{\delta}$ contains $p_0$.  We prove in this case that the $(\widetilde{G}, \tEin)$-structure on $M$ is complete---that is, $\delta$ is a diffeomorphism.  

\begin{theorem}
\label{thm.case1final}
If the image of $\bar{\delta}$ contains $p_0$, then up to a finite covering, $M$ is the quotient of $\tEin$ by a free, properly discontinuous $\BZ$-action that leaves $\Delta$ invariant.  More precisely, $M$ is finitely covered by $\tEin \backslash \langle g \alpha^i \rangle$ where $i > 0$ and $g \in \widetilde{\mathcal{U}}^0$.
\end{theorem}

The first step of the proof deals with the developing map \emph{vis-a-vis} $\Delta \subset \tEin$; recall $\Delta$ is connected.  The following proposition serves for this case as well as the second case in the next section.  Thus there is no special assumption on $\delta$ at this stage.
\begin{proposition}\label{prop:diffonLambda}
Any connected component $\Lambda$ of $\delta^{-1}(\Delta)$ is mapped by $\delta$ diffeomorphically onto its image. In particular,
\begin{enumerate}
   \item if $p_0 \in \bar{\delta}(\Lambda)$, then $\delta$ maps $\Lambda$ diffeomorphically onto $\Delta$.
    \item if $p_0 \notin \bar{\delta}(\Lambda)$, then $\delta$ maps $\Lambda$ diffeomorphically onto a connected component of $\Delta \setminus \{ \tilde{p}_i \}$, where $ \{ \tilde{p}_i \} = \pi_{\mathrm{Ein}}^{-1}(p_0)$.
\end{enumerate}
\end{proposition}

\begin{proof}
The restriction of $\mathcal{U}$ to $\bar{\Delta}$ is the parabolic flow with unique fixed point $p_0$---see \S \ref{sec.unip.subgp}.  Denote the corresponding vector field on $\bar{\Delta}$, and its lift to $\Delta$, by $Y_\sigma$. The lift to $\Delta$ vanishes precisely at the points $\{ \tilde{p}_i \}$ and has no periodic points in each connected component of $\Delta \backslash \{ \tilde{p}_i \}$.

By Proposition \ref{prop.vf.pullback}, $\delta^* Y_\sigma$ restricted to $\delta^{-1}(\Delta)$ is a complete vector field on $\delta^{-1}(\Delta)$, and on $\Lambda$ in particular; moreover, the image $\delta(\Lambda)$ is invariant by the flow along $Y_\sigma$.  The flow-invariant subsets of $\Delta$ are subsets of $\{ \tilde{p}_i \}$ and unions of components of $\Delta \backslash \{ \tilde{p}_i \}$.  

In case (2) of this proposition, $\delta(\Lambda)$ is one component of $\Delta \backslash \{ \tilde{p}_i \}$, which we will call $\Delta_i$.  
Such a component is parametrized by the flow along $Y_\sigma$ of any point in it.  As $\delta$ intertwines the two corresponding flows, it follows in this case that $\delta$ maps $\Lambda$ diffeomorphically onto $\Delta_i$, as desired.

In case (1), observe that the saturation of $\Lambda$ by the deck transformations of $\widetilde{M} \rightarrow M$ is a union of connected components of the closed set $\delta^{-1}(\Lambda)$, and therefore is also closed.  Then $\delta(\Lambda) \subset M$ is a closed photon, which inherits from the $(G, \Ein)$-structure on $M$ a $(\mbox{PSL}(2,\BR), {\bf RP}^1)$-structure.  Of course, the latter structure also has unipotent holonomy.  
From the classification of one-dimensional $(\mathrm{PSL}(2,\BR)$, $\RP^1)$-manifolds---see \cite{kuiper.RP1} or \cite[\S 5.5]{goldman.gstom}---the structure on $\delta(\Lambda)$ is complete.  The developing map of this structure is a diffeomorphism to $\widetilde{\RP^1}$, which is embedded into $\tEin$ as $\Delta$---see \S \ref{sec.geom.tein}. This developing map is the restriction of $\delta$ to $\Lambda$, so we conclude that $\delta$ maps $\Lambda$ diffeomorphically onto $\Delta$.
\end{proof}

Completeness will extend from $\Lambda$ to all of $\widetilde{M}$ with the help of the $\tau$-flow.  By Proposition \ref{prop.vf.pullback}, $\delta^* Y_\tau = X_\tau$ is a complete vector field on $\widetilde{M}$ and the image of $\delta$ is $\tau$-saturated.  Denote by $\{ \hat{\tau}^s \}$ the corresponding flow on $\widetilde{M}$.  The following lemma is proved using only the equivariance and local diffeomorphism properties of $\delta$.  The interval $I$ in the statement need not map diffeomorphically onto its image.

\begin{lemma}\label{Mlimitopen}
Let $I\subseteq\delta^{-1}(\Delta)$ be open and connected. Then the set
\begin{equation*}
    W^+(I) = \{x\in\widetilde{M}\setminus\delta^{-1}(\Delta): \lim_{t\to\infty}\hat{\tau}^t. x \in I\}
\end{equation*}
is open in $\widetilde{M}$.  The same holds for the analogously defined set $W^-(I)$.
\end{lemma}

\begin{proof}
Let $x \in W^+(I)$ with $x_\infty= \lim_{t\to\infty}\hat{\tau}^t.x$, and let $y = \delta(x)$.  The developing map intertwines the $\hat{\tau}$- and $\tau$-flows.  Because $\lim_{t\to\infty}\hat{\tau}^t.x$ exists, it is mapped under $\delta$ to $\lim_{t \to \infty} \tau^t.y$ by continuity.   By assumption, $y \notin \Delta$, so this latter limit can be expressed as $\rho_\Delta(y)$.  Choose connected neighborhoods $x_\infty \in A \subset \widetilde{M}$ and $\rho_\Delta(y) \in B \subset \tEin$ such that $\delta$ maps $A$ diffeomorphically to $B$; shrink them if necessary to ensure that $A \cap \Lambda \subset I$.

Choose connected neighborhoods $x \in U \subset \widetilde{M}$ and $y \in V \subset \tEin \backslash \Delta$ such that $\delta$ maps $U$ diffeomorphically to $V$, and let
$$ U' = \left( \rho_\Delta \circ \left. \delta \right |_U \right)^{-1}(B \cap \Delta)$$
which is again an open neighborhood of $x$.  Let $V' = \delta(U')$.

Consider $x' \in U'$, and let $y' = \delta(x')$.  As $\{ \tau^t \}$ converges uniformly on compact sets to $\rho_\Delta$, there are $T > 0$ and a connected open $W \subset \overline{W} \subset V'$ containing $y$ and $y'$ with $\tau^t(W) \subset B$ for all $t \geq T$.  We can choose $T$ large enough that $\hat{\tau}^T(x) \in A$.  Connectedness of $\{ \tau^t.y : t \geq T\} \subset B$ implies that $\hat{\tau}^t.x \in A$ for all $t \geq T$.  Now $(\left. \delta \right|_{U'})^{-1}(W)$ is connected and open, so for any $t \geq T$,
$$\hat{\tau}^t((\left. \delta \right|_{U'})^{-1}(W)) = ( \left. \delta \right|_A)^{-1}(\tau^t(W)) \subset A$$
Since $\lim_{t \rightarrow \infty} \tau^t(y')$ exists and belongs to $B \cap \Delta$, it now follows via equivariance that $\lim_{t \rightarrow \infty} \hat{\tau}^t(x')$ exists and belongs to $A \cap \Lambda  \subset I$.  As $x'$ was an arbitrary point of $U'$, we conclude that $U' \subset W^+(I)$, which completes the proof for $W^+(I)$.  The proof for $W^-(I)$ is completely analogous.  
%
\end{proof}

Adding the assumption that $I$ maps diffeomorphically onto its image brings us to the key step for proving completeness.

\begin{proposition}\label{limtolim}
Let $J\subseteq\Delta$ be open and connected and $I$ be a connected component of $\delta^{-1}(J)$. If $\delta$ maps $I$ to $J$ diffeomorphically, then the sets
\begin{equation*}
    W^\pm(I) = \{x\in\widetilde{M}\setminus\delta^{-1}(\Delta): \lim_{t\to\pm\infty}\hat{\tau}^t.x \in I\}
\end{equation*}
each map diffeomorphically onto the sets
\begin{equation*}
    \Omega^\pm(J) = \{z\in\tEin\setminus\Delta: \lim_{t\to\pm\infty}\tau^t.z \in J\}.
\end{equation*}
\end{proposition}

\begin{proof}
The proof here is for $W^+(I)$ and $\Omega^+(J)$; the other case is proved \emph{mutatis mutandis}.

Let $x_1, x_2 \in W^+(I)$ with $y_i = \delta(x_i)$, $i=1,2$.   As $\lim_{t \rightarrow \infty} \hat{\tau}^t .x_i = z_i$ exists, it maps under $\delta$ to $\lim_{t \rightarrow \infty} \tau^t.y_i$, for $i=1,2$.  
Suppose that $y_1 = y_2 = y$.  Then by equivariance
$$\delta(z_1) = \delta(z_2) = \lim_{t \rightarrow \infty} \tau^t(y)$$
By our assumption, this implies $z_1 = z_2 = z$.  Now let $A$ be a neighborhood of $z$ mapping diffeomorphically under $\delta$ to its image, which we will denote $B$. Let $T > 0$ be such that $\hat{\tau}^T(x_1), \hat{\tau}^T(x_2) \in A$.  By equivariance, $\tau^T(y) \in B$.  But then $\hat{\tau}^T(x_1)= \hat{\tau}^T(x_2)$ which implies $x_1 = x_2$.  Therefore $\delta$ is injective in restriction to $W^+(I)$.

Next let $y \in \Omega^+(J)$, with $\lim_{t\to\infty}\tau^t.y=z$. Let $w \in I$ be the $\delta$-preimage of $z$.  Let $A$ be a neighborhood of $w$ mapping diffeomorphically to $B$, a neighborhood of $z$.  Let $T > 0$ such that $\tau^t.y \in B$ for all $t\geq T$. Let $x = (\left. \delta \right|_A)^{-1}(\tau^T.y)$.  Then $\lim_{t \rightarrow \infty} \hat{\tau}^t(x) = w$ and $\delta(\tau^{-T}.x) = y$.

Now $\delta$ restricted to the open set $W^+(I)$ is a bijective local diffeomorphism onto $\Omega^+(J)$, hence a diffeomorphism onto $\Omega^+(J)$.
\end{proof}

We are ready to assemble the proof of Theorem \ref{thm.case1final}.  Assume that $p_0$ is in the image of $\bar{\delta}$.  Let $\Lambda$ be a connected component of $\delta^{-1}(\Delta)$.  By Proposition \ref{prop:diffonLambda} (2), $\delta$ maps $\Lambda$ diffeomorphically onto $\Delta$.  Let $\Omega = \tEin \backslash \Delta$ and $W = W^+(\Lambda)$.  Proposition \ref{limtolim} above says that $W$ maps diffeomorphically under $\delta$ to $\Omega$.  Since $\dim M \geq 3$, the codimension of $\delta^{-1}(\Delta)$ is at least two, so $\widetilde{M} \backslash \delta^{-1}(\Delta)$ is path connected.  
%
%
Lemma \ref{daly4} applies to $U = \widetilde{M} \backslash \delta^{-1}(\Delta)$ and $X = \Omega$, to give that $W = W^+(\Lambda)$ maps diffeomorphically to $\Ein\setminus \Delta$ and equals $U$. Any other component $\Lambda'$ of $\delta^{-1}(\Delta)$ would give $W^+(\Lambda')$ disjoint from $W^+(\Lambda)$ yet equal to $U$---a contradiction.  Therefore $\Lambda = \delta^{-1}(\Delta)$, so $\widetilde{M} = W \cup \Lambda$.
 Finally, $\delta$ is a bijective local diffeomorphism, so $\delta$ maps $\widetilde{M}$ diffeomorphically onto $\tEin$.  Completeness is proved.

To analyze the holonomy, recall that $\pi_M(\Lambda) \cong \bar{\Delta}$ is a complete $(\mbox{PSL}(2,\BR), {\bf RP}^1)$-manifold.  The restricton to $\Lambda$ of the covering group action on $\widetilde{M}$ is conjugated by $\delta$ to the action of $\langle \alpha^i \rangle$ on $\Delta$, for some $i > 0$.   On the other hand, this restriction to $\Lambda$ is faithful, so $\pi_1(M) \cong\BZ$.  Since we assume that the holonomy projected to $\PO(n,2)$ belongs to $\mathcal{U}$, it follows that 
$\Gamma = \langle \alpha^i g \rangle$ with $i > 0$ and $g \in \widetilde{\mathcal{U}}^0$.

\section{The developing image does not contain $p_0$ but does meet $\bar{\Delta}$}
\label{sec.case2}

We continue to the second case of our classification. 
This case does not have an analogue in the conformal Riemannian setting. 

\subsection{The developing image}

Write $\{ \tilde{p}_i \} = \pi_{\mathrm{Ein}}^{-1}(p_0)$, and let $\Delta_0$ be a connected component of $\Delta \backslash \{ \tilde{p}_i \}$ meeting the image of $\delta$.  Let $\Lambda_0$ be a connected component of $\delta^{-1}(\Delta_0)$.  By Proposition \ref{prop:diffonLambda} (2), $\Lambda_0$ maps diffeomorphically under $\delta$ onto $\Delta_0$.  The latter set is of the form 
$$ \Delta_0 = \{ x \in \Delta \ : \ \tilde{p}_i < x < \alpha. \tilde{p}_i \}$$
for some $i$, which we will assume to be $0$.   

Let $\hat{\Gamma}_0 < \Conf \widetilde{M}$ be the stabilizer of $\Lambda_0$ in the group of deck transformations of $\widetilde{M} \rightarrow M$.  As $\Delta_0$ is homeomorphic to $\BR$, so is $\Lambda_0$.  Proposition \ref{prop.closed.invt} gives that $\pi_M(\delta^{-1}(\Delta))$ is closed; then $\pi_M(\Lambda_0)$, a connected component of this set, is also closed. 
Thus $\hat{\Gamma}_0 \cong \BZ$.

Let $S_0 \subset L(\tilde{p}_0) \backslash \{ \tilde{p}_i \}$ be the connected component containing $\Delta_0$.  Note that $S_0 = L(\tilde{p}_0,\alpha. \tilde{p}_0)$.  Let $\Sigma_0$ be a connected component of $\delta^{-1}(S_0)$.  The following proposition establishes completeness of $\delta$ between $\Sigma_0$ and $S_0$.  It is not assumed that $\delta(\Sigma_0)$ meets $\Delta$; part (2) of the proposition will be used for case 3, in \S \ref{sec.case3} below.

\begin{proposition}
  \label{prop.lightcone.complete}
  Assume that $p_0 \notin \mbox{im } \bar{\delta}$.
  Let $S$ be a connected component of $L(\tilde{p}_0) \backslash \pi_{\mathrm{Ein}}^{-1}(p_0)$ and $\Delta_0 = S \cap \Delta$.  Let
  $\Sigma$ be a connected component of $\delta^{-1}(S)$.  Then
  $\delta$ maps $\Sigma$ diffeomorphically to its image, which equals
    \begin{enumerate}
    \item $S$ if $\delta(\Sigma) \cap \Delta_0 \neq \emptyset$.
       \item $S \backslash \Delta_0$ if $\delta(\Sigma) \cap
      \Delta_0 = \emptyset$
             \end{enumerate}
  \end{proposition}

  \begin{proof}
  The submanifold $\Sigma$ is a connected component of $\delta^{-1}(L(\tilde{p}_0))$, so it is closed in $\widetilde{M}$.  On the other hand, $S$ is a regular submanifold.  
The holonomy subgroup $\Gamma_0 = h(\hat{\Gamma}_0)$ leaves $S$ invariant, as does all of $\widetilde{\mathcal{U}}^0$.  The $\Gamma$-orbit of $S$ is a union of connected components of $L(\tilde{p}_0) \backslash \{ \tilde{p}_i \}$.  Thus $\Sigma$ can also be considered a connected component of the inverse image of the $\Gamma$-invariant regular submanifold $\Gamma.S$.  

Let $\mathcal{F}$ be the foliation of $S$ by photons.  It is invariant by $\Gamma_0$ (and extends to a $\Gamma$-invariant foliation of $\Gamma.S$).  
We have seen just above that if $\delta(\Sigma)$ meets $\Delta_0$ then any connected component of $\delta^{-1}(\Delta_0)$ in $\Sigma$ maps diffeomorphically to $\Delta_0$ under $\delta$.  On any photon in $S\backslash \Delta_0$, the $\tau$-flow acts simply transitively---this can be seen from the formula (\ref{eqn.tau.flow}) for the $\tau$-flow, taking $x_{n+1} = 0$ and $x_n \neq 0$.  By Proposition \ref{prop.vf.pullback}, $\delta^* Y_\tau = X_\tau$ is a complete vector field on $\widetilde{M}$.  It follows that for any photon $\gamma$ of $S\backslash \Delta_0$, each component of $\delta^{-1}(\gamma) \cap \Sigma$ is mapped diffeomorphically by $\delta$ to $\gamma$.  We conclude that $\delta$ maps leaves of $\hat{\mathcal{F}}$ in $\Sigma$ diffeomorphically to leaves of $\mathcal{F}$ in $S$, where $\hat{\mathcal{F}}$ is the pulled-back foliation by photons on $\Sigma$.

Next we define vector fields $Y_2, \ldots, Y_{n-1}$ on $S$.  They will be lifted from $\bar{S} = L(p_0) \backslash \{ p_0 \} \subset \Ein$.  Denote by $\bar{\mathcal{F}}$ the foliation by photons on $\bar{S}$.
Consider the map given in homogeneous coordinates by
\begin{equation}
\label{eqn.param.lightcone}
 \iota : \BR \times \BR^{n-2} \rightarrow \bar{S} \qquad (t,(y_2, \ldots, y_{n-1})) \mapsto [t: - \frac{||y||^2}{2} : y_2 : \cdots : y_{n-1}: 1: 0]
 \end{equation}
 It is an injective immersion onto $\bar{S} \backslash \bar{\Delta}$.  Let $Y_i = \iota_*(\partial y_i)$ for $i=2, \ldots, n-1$.  As $\partial y_i$ is projectable under $\BR \times \BR^{n-2} \rightarrow \BR^{n-2}$, and $\iota$ maps the fibers $\BR \times \{y \}$ diffeomorphically to the photons of $\bar{S}\backslash \bar{\Delta}$, it follows that $Y_i$ is projectable modulo $\bar{\mathcal{F}}$ for all $i$.  The $\mathcal{U}$-action on $\Ein$ preserves the image of $\iota$ and the foliation $\mathcal{F}$.  The conjugated $\mathcal{U}$-action on $\BR \times \BR^{n-2}$ descends to the $\BR^{n-2}$-action by translations on $(\BR \times \BR^{n-2})/\BR \cong \BR^{n-2}$, and thus centralizes $\partial y_i$ modulo the $\BR$-factor, for all $i$.  The $\mathcal{U}$-action on $\bar{S} \backslash \bar{\Delta}$ thus centralizes each $Y_i$ modulo $T \bar{\mathcal{F}}$.

Let $L = L([e_0]) \cap L([e_{n+1}])$, a conformally embedded copy of ${\bf S}^{n-2}$, corresponding in homogeneous coordinates on $\Ein$ to the locus where $x_0=0=x_{n+1}$.  In the affine chart on ${\bf P}(e_0^\perp \cap e_{n+1}^\perp) \cong {\bf RP}^{n-1}$ where $x_n=1$, a routine calculation using (\ref{eqn.param.lightcone}) gives $Y_i = \partial x_i - x_i \partial x_1$.  In the affine chart where $x_1=1$,  the expression is
$$Y_i = \sum_{j=2}^{n-1} x_i x_j \partial x_j + x_n \partial x_i + x_n x_i \partial x_n$$
which tends to $0$ as $x \rightarrow 0$.  The origin of this affine chart corresponds to $[e_1] \in \bar{\Delta}$.  The conclusion is that $Y_i$ extends by $0$ to $\bar{\Delta}$.  Now lift the resulting vector fields on $\bar{S}$ to $S$; we will continue to denote these by $Y_i$, for $i=2, \ldots, n-1$.

The vector fields $Y_i $ are each projectable and $\widetilde{U}^0$-invariant modulo $\mathcal{F}$, thus they are $\Gamma_0$-invariant (and $\Gamma$-invariant when equivariantly extended to $\Gamma.S$) modulo $\mathcal{F}$.  We are now in a position to apply Proposition \ref{prop.fgh.technique}. The image $\delta(\Sigma)$ is invariant by the flow along $Y_i$ for all $i$, which means that it meets every leaf of $S \backslash \Delta_0$.  In case (1) or (2), $\delta$ maps $\Sigma$ onto the claimed set.  In case (2), the vector fields $\{Y_i \}$ form a framing of the local leaf space at every point of $S \backslash \Delta_0$, so part (2) of Proposition \ref{prop.fgh.technique} applies to give that $\Sigma \rightarrow S \backslash \Delta_0 \cong \BR^{n-1}$ is a covering map, which means it is a diffeomorphism.

In case (1), we take $W = S \backslash \Delta$ in Proposition \ref{prop.fgh.technique} (2).  It gives that any connected component of $\Sigma \backslash \delta^{-1}(\Delta_0)$ maps by a covering map onto $S \backslash \Delta_0 \cong \BR^{n-1}$, therefore by a diffeomorphism.   By the following lemma, there is only one such component.  Then $\delta$ is a bijective local diffeomorphism on $\Lambda \cup \Sigma \backslash \delta^{-1}(\Delta_0) = \Sigma$, so it maps $\Sigma$ diffeomorphically to $S$, as desired.
    \end{proof}

\begin{lemma}
For $S$ and $\Sigma$ as in Proposition \ref{prop.lightcone.complete} case (1),
 $\Sigma \backslash \delta^{-1}(\Delta_0)$ is connected.
\end{lemma}

\begin{proof}
Let $\Lambda$ be a connected component of $\delta^{-1}(\Delta_0) \cap \Sigma$.  Consider $W^\pm(\Lambda)$ as in Lemma \ref{Mlimitopen}.  Since $\Lambda$ maps diffeomorphically under $\delta$ to $\Delta_0$,  Proposition \ref{limtolim} gives that $W^\pm(\Lambda)$ maps diffeomorphically to $\Omega^\pm(\Delta_0)$, respectively.  By Proposition \ref{prop.minlimit}, these are $\mbox{\bf Min}^{\mp}(\tilde{p}_0)$, respectively.  By Corollary \ref{Minboundary}, the boundaries of these sets in $\tEin \backslash \{ \tilde{p}_i \}$ are $\alpha^{-1}.S \cup S$ and $S \cup \alpha.S$, respectively.

For $x \in \Sigma \cap \partial W^+(\Lambda)$, let $\hat{U}$ be a neighborhood mapping diffeomorphically to its image under $\delta$, which will be denoted $U$.  Now
$$ \hat{U} \cap \partial W^+(\Lambda) = \hat{U} \cap \delta^{-1}(\partial \Omega^+(\Delta_0) ) = \hat{U} \cap \delta^{-1}(\alpha^{-1}.S \cup S)$$
 By shrinking $\hat{U}$ to ensure that $U \cap (L(\tilde{p}_0) \backslash \{ \tilde{p}_i \})$ is connected, we may arrange that 
 $$ \hat{U} \cap \partial W^+(\Lambda) = \hat{U} \cap  \delta^{-1}(S) = \hat{U} \cap \Sigma$$
 Thus $\Sigma \cap \partial W^+(\Lambda)$ is open in $\Sigma$.  Since it is also closed and $\Sigma$ is connected, it follows that $\Sigma \subset \partial W^+(\Lambda)$.  The analogous argument implies $\Sigma \subset \partial W^-(\Lambda)$.  

Now suppose $x \in \Sigma \backslash \delta^{-1}(\Delta_0)$.  Then $\hat{U}$ meets $W^\pm(\Lambda)$ and $U$ meets only the Minkowski patches $\mbox{\bf Min}^\pm(\tilde{p}_0)$; by Proposition \ref{limtolim}, $\hat{U} \backslash \Sigma \subset W^+(\Lambda) \cup W^-(\Lambda)$.  Now $\delta$ is a local diffeomorphism between the open sets
 $$W^+(\Lambda) \cup (\Sigma \backslash \delta^{-1}(\Delta_0) ) \cup W^-(\Lambda) \rightarrow \mbox{\bf Min}^-(\tilde{p}_0) \cup (S \backslash \Delta_0) \cup \mbox{\bf Min}^+(\tilde{p}_0)$$
Any component $\Sigma'$ of $\Sigma \backslash \delta^{-1}(\Delta_0)$ maps diffeomorphically to $S \backslash \Delta_0$, as seen in the proof of Proposition \ref{prop.lightcone.complete}.  Then the open subset $W^+(\Lambda) \cup \Sigma' \cup W^-(\Lambda)$ maps diffeomorphically onto the target above.  Therefore $\Sigma \backslash \delta^{-1}(\Delta_0) = \Sigma'$ by Proposition \ref{daly4}.
\end{proof}

With Proposition \ref{prop.lightcone.complete} and the arguments of the preceding proof, we have also established:
\begin{proposition}
\label{prop.W.Omega}
Let $S$ and $\Sigma$ be as in case (1) of Proposition \ref{prop.lightcone.complete}.  Then 
\begin{itemize}
\item   $W = W^+(\Lambda) \cup \Sigma \cup W^-(\Lambda)$ is open and is mapped diffeomorphically by $\delta$ onto $\Omega = \mbox{\bf Min}^-(\tilde{p}_0) \cup S \cup \mbox{\bf Min}^+(\tilde{p}_0)$.
\item $\delta^{-1}(\Delta_0) \cap \Sigma = \Lambda$ is connected.
\end{itemize}
\end{proposition}

The following complete description of $\delta$ will conclude this subsection.

\begin{proposition}
The set $W$ in the conclusion of Proposition \ref{prop.W.Omega} equals $\widetilde{M}$, and $\Gamma_0 = \Gamma$. The developing map is a diffeomorphism of $\widetilde{M}$ onto $\Omega$, and the holonomy image is a subgroup of $\widetilde{\mathcal{U}}^0$ isomorphic to $\BZ$.
\end{proposition}

\begin{figure}
  \centering
  \includegraphics[]{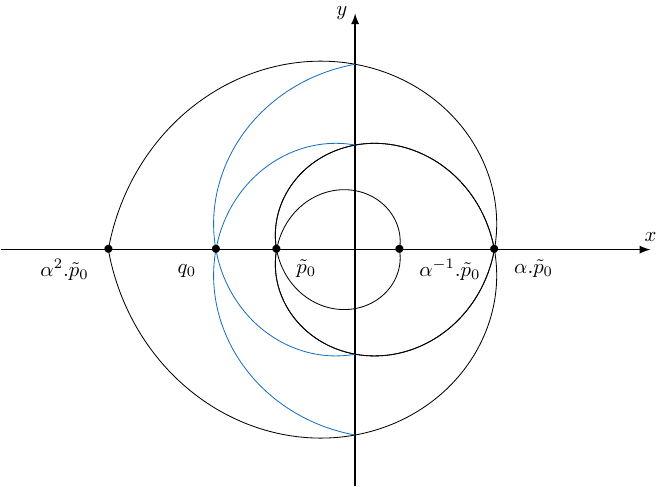}
  \caption{The lightcone $L(q_0)$ in $\mbox{\bf Min}^+(\tilde{p}_0)$ has compact closure in $\tEin \backslash \{ \tilde{p}_i \}$.}
  \label{fig:afflightcones}
\end{figure}

\begin{proof}
Assume to the contrary that there is $x \in \partial W$. Then $\delta(x) \in \partial \Omega \cap (\tEin\setminus \{ \tilde{p}_i \})$.  From Corollary \ref{Minboundary} the boundaries of $\mbox{\bf Min}^\pm(\tilde{p}_0) \cap \tEin\setminus \{ \tilde{p}_i \}$ are contained in the three punctured light cone components $L(\alpha^i.\tilde{p}_0, \alpha^{i+1}.\tilde{p}_0) $ for $i=-1, 0, 1$.  One of these, $L(\tilde{p}_0,\alpha.\tilde{p}_0) $, equals $S$ and is in the interior of $\Omega$.  Therefore 
$$ \partial \Omega = L(\alpha^{-1}.\tilde{p}_0,\tilde{p}_0) \cup L(\alpha.\tilde{p}_0, \alpha^2.\tilde{p}_0) = \alpha^{-1}.S \cup \alpha.S$$

By construction, $\Gamma_0$ is contained in the stabilizer of $\Delta_0$; the latter intersects $\widetilde{\mathcal{U}}$ in $\widetilde{\mathcal{U}}^0$.  Thus $\Gamma_0$ is also contained in the stabilizer of $\alpha^i.S$ and $\alpha^i.\Delta_0$ for all $i$.
Without loss of generality, we assume $\delta(x) \in \alpha.S$. Let $\Sigma' \subset \delta^{-1}(\alpha.S)$ be a nonempty connected component. Proposition \ref{prop.lightcone.complete} implies that $\Sigma'$ maps diffeomorphically onto $\alpha.S$ or $\alpha.S \backslash \alpha.\Delta_0$ under $\delta$.  By the usual argument with Proposition \ref{prop.closed.invt}, the stabilizer in $\pi_1(M)$ of $\Sigma'$ acts cocompactly.  On the other hand, this stabilizer maps isomorphically to $\Gamma_0 \cong \BZ$.  It is then impossible that $\delta(\Sigma') = \alpha.S \backslash \alpha.\Delta_0 \cong \BR^{n-1}$.  Therefore necessarily $\Sigma'$ maps diffeomorphically under $\delta$ onto $\alpha.S$.

Now $\Sigma' \cup W^-(\Lambda) \cup \Sigma$ develops diffeomorphically to $\alpha.S \cup \mbox{\bf Min}^+(\tilde{p}_0) \cup S$. The latter set is the closure of $\mbox{\bf Min}^+(\tilde{p}_0)$ in $\tEin\setminus\{ \tilde{p}_i \}$ by Corollary \ref{Minboundary}.  We will show that $\Gamma_0 \cong \BZ$ does not act properly on this set.
Let $q_0$ be the origin of the Minkowski patch $\mbox{\bf Min}^+(\tilde{p}_0)$, and let $K = L(q_0) \cap \mbox{\bf Min}^+(\tilde{p}_0)$. Since any two lightcones in $\Min$ intersect, the images $\gamma.K$ intersect $K$ for any  $\gamma \in \Gamma_0$. The closure 
$\overline{K} \subset \overline{\mbox{\bf Min}^+(\tilde{p}_0)}$ does not meet $\{ \tilde{p}_i \}$---it is easily seen from the Minkowski chart in \S \ref{sec.geom.ein} that $p_0 \notin \pi_{\mathrm{Ein}}(\overline{K})$---so $\overline{K}$ is a compact subset of $\tEin\setminus\{ \tilde{p}_i \}$ contained in $\alpha.S \cup \mbox{\bf Min}^+(\tilde{p}_0) \cup S$.
Thus $\delta^{-1}(K) \cap \left( \Sigma' \cup W^-(\Lambda) \cup \Sigma \right)$ is compact, and intersects its image under any $\gamma \in \hat{\Gamma}_0 \cong \BZ$.
Because $\hat{\Gamma}_0$ acts properly by deck transformations on $\widetilde{M}$, this is a contradiction.
See Figure \ref{fig:afflightcones}.   The conclusions of the proposition now follow.
\end{proof}

\subsection{Determination of holonomy, conclusion of classification}

It is established that $\hat{\Gamma}_0 \cong \Gamma \cong \BZ$; moreover, $\Gamma$ is generated by an element of $\widetilde{\mathcal{U}}^0$ acting freely, properly discontinuously, and cocompactly on $\Delta_0$, $S$, and $\Omega$.  
We first establish necessary conditions on this generator, which we will call $\gamma$.

\begin{proposition}
\label{prop.case2.holo}
Let $\gamma$ be the generator of the holonomy image under the assumptions of this section, and let $\bar{\gamma} = q(\gamma) \in G$.  Let $L(\bar{\gamma}) + u_{\bar{\gamma}}$ be the affine decomposition as in Section \ref{sec.unip.subgp}.  Then $u_{\bar{\gamma}}$ is nontrivial modulo $e_1^\perp$.  Moreover, if $L(\bar{\gamma}) = \mbox{Id}$, then $u_{\bar{\gamma}}$ is timelike.
\end{proposition}  

\begin{proof}
The restriction of $\mathcal{U}$ to $\bar{\Delta} = {\bf P}(\mbox{span} \{ e_0, e_1 \})$ corresponds in the 
affine representation to the projection of $u_{\bar{\gamma}}$ modulo $e_1^\perp$, and this must be nontrivial.

Next suppose that $\bar{\gamma} \in \ker L$, and let $v = u_{\bar{\gamma}}$.   Let $\iota : \BR \times \BR^{n-2} \rightarrow \bar{S}'$ be the chart on $\bar{S}' = \pi_{\mathrm{Ein}}(S)\backslash \bar{\Delta}$ given in (\ref{eqn.param.lightcone}).  In this chart, the action of $\bar{\gamma}$ is by 
$$(t,y) \mapsto (t + \langle v, (- ||y||^2/2, y_2, \ldots, y_n, 1) \rangle, y)$$
where the scalar product on $\BR^{n-1,1}$ is the one with quadratic form $q_{n-1,1}$.  The lines 
$$ \{ s  (- ||y||^2/2, y_2, \ldots, y_n, 1) \ : s \in \BR, \ y \in \BR^{n-2} \}$$ 
describe the full null cone $\mathcal{N}^{n-1,1}$ of $\BR^{n-1,1}$ except for the line $\BR e_1$.  We have already established that $v$ is not orthogonal to $e_1$.  Then $\bar{\gamma}$ will pointwise fix some photon of $\bar{S}$ unless $v$ is not orthogonal to any line in $\mathcal{N}^{n-1,1}$.  Since $\bar{\gamma}$ must act freely on $\bar{S}$, we conclude that $v$ is timelike.
\end{proof}

\begin{lemma}\label{fundconds}
Let $q_0 \in \mbox{\bf Min}^-(\tilde{p}_0)$, and suppose there is $\gamma \in \widetilde{G}$ that satisfies
\begin{enumerate}
    \item $q_0 \ll \gamma. q_0$,
    \item $\displaystyle \lim_{i\to\infty}\gamma^i .q_0 =\alpha.\tilde{p}_0$ and $\displaystyle \lim_{i\to-\infty}\gamma^i. q_0= \alpha^{-1}.\tilde{p}_0$.
\end{enumerate}
Then $\langle \gamma \rangle$ acts on $\Omega$ with precompact fundamental domain given by
\begin{equation*}
    D = J^+(q_0)\cap I^-(\alpha \gamma .q_0).
\end{equation*}
\end{lemma}

\begin{figure}
  \centering
  \includegraphics[]{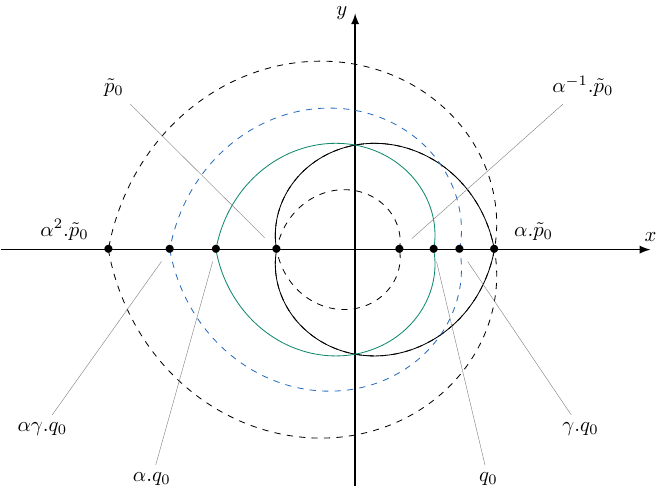}
  \caption{The fundamental domain for the action of $\langle \gamma \rangle$ as in Lemma \ref{fundconds}.}
  \label{fig:funddomain}
\end{figure}

See Figure \ref{fig:funddomain} for a diagram of this fundamental domain in $\tEin$.

\begin{proof} We first prove that the union $\displaystyle \bigcup_{i\in\BZ}\gamma^i .D$ covers the entirety of $\Omega$. From Lemma \ref{Eincomplement}, 
\begin{align*}
    \bigcup_{i\in\BZ} \gamma^i .(J^+(q_0)\cap I^-(\alpha \gamma .q_0))
    &= \bigcup_{i\in\BZ} \gamma^i .(I^-(\alpha \gamma. q_0)\setminus I^-(\alpha .q_0))\\
    &= \bigcup_{i\in\BZ} (I^-(\alpha \gamma^{i+1} .q_0)\setminus I^-(\alpha \gamma^i .q_0))
\end{align*}
As both $\alpha$ and $\gamma$ preserve time orientation, 
\begin{equation*}
    \alpha \gamma^i .q_0 \ll \alpha \gamma^{i+1} .q_0,\quad\forall i\in\BZ
\end{equation*}
which gives
\begin{equation*}
    I^-(\alpha \gamma^i. q_0) \subset I^-(\alpha \gamma^{i+1} .q_0),\quad\forall i\in\BZ.
\end{equation*}
It thus follows that
\begin{align*}
    \bigcup_{i\in\BZ} (I^-(\alpha \gamma^{i+1}. q_0)\setminus I^-(\alpha \gamma^i .q_0))
    &= \bigcup_{i\in\BZ}I^-(\alpha \gamma^i .q_0)\setminus \bigcap_{j\in\BZ} I^-(\alpha\gamma^j .q_0)
\end{align*}
 From (2),
\begin{center}
    $\displaystyle \lim_{i\to\infty}\alpha\gamma^i .q_0 =\alpha^2.\tilde{p}_0$.
\end{center}
Since $\alpha^2 .\tilde{p}_0$ succeeds each $\alpha\gamma^i .q_0$, 
\begin{equation*}
   \bigcup_{i \in \BZ} I^-(\alpha\gamma^i .q_0)\subseteq I^-(\alpha^2 .\tilde{p}_0) 
\end{equation*}
For any set $S$,  the future and past satisfy $I^\pm[\overline{S}] = I^\pm[S]$ \cite[Prop 2.11]{penrose.diff.top.rel}.  Taking $S = \{ \alpha \gamma^i.q_0 \} _{i \in \BZ }$, we obtain the reverse containment and conclude
\begin{equation*}
    \bigcup_{i\in\BZ}I^-(\alpha \gamma^i .q_0) = I^-(\alpha^2 .\tilde{p}_0).
\end{equation*}

From Lemma \ref{Eincomplement},
\begin{equation*}
    \big(\bigcap_{j\in\BZ}I^-(\alpha\gamma^j .q_0)\big)^c = \bigcup_{j\in\BZ} J^+(\gamma^j .q_0).
\end{equation*}
For any $x > \gamma^i .q_0$, by (1),  $\gamma^{i-1}.q_0 \ll x$ and $J^+(\gamma^i .q_0)\subset I^+(\gamma^{i-1}.q_0)$.  It follows that
\begin{align*}
    \bigcup_{j\in\BZ} J^+(\gamma^j .q_0) 
    &= \bigcup_{j\in\BZ} I^+(\gamma^j .q_0)
  \end{align*}

From (2), $\displaystyle \lim_{i\to-\infty}\gamma^i .q_0= \alpha^{-1}.\tilde{p}_0$.  For all $i$, $\gamma^i .q_0 \in \mbox{Min}^-(\tilde{p}_0)\subset I^+(\alpha^{-1}.\tilde{p}_0)$, so $\gamma^i.q_0  > \alpha^{-1}.\tilde{p}_\infty$ for all $i$.  By the argument with \cite[Prop 2.11]{penrose.diff.top.rel} as above, we obtain
$$\bigcup_{j\in\BZ} I^+(\gamma^j .q_0) = I^+(\alpha^{-1}.\tilde{p}_\infty)$$
We can finally conclude that 
\begin{equation*}
\bigcup_{i\in\BZ}\gamma^i .D = I^-(\alpha^2.\tilde{p}_0) \cap I^+( \alpha^{-1}.\tilde{p}_0)
\end{equation*}
which by definition equals $\mbox{\bf Min}^-(\tilde{p}_0) \cup \mbox{\bf Min}^+(\tilde{p}_0) = \Omega$.

It remains to show $D\cap \gamma^i .D = \emptyset$ for all $i \in \BZ$.  Let $i$ be given, which we may assume is positive.  From (1), it follows that
\begin{equation*} 
    J^+(\gamma^i .q_0) \subset J^+(\gamma. q_0).
\end{equation*}
By Lemma \ref{Eincomplement}, we can also express
\begin{equation*}
    D = J^+(q_0) \setminus J^+(\gamma .q_0) \quad \mbox{and} \quad 
    \gamma^i. D = J^+(\gamma^i. q_0)\setminus J^+(\gamma^{i+1}.q_0).
\end{equation*}
The desired disjointness now follows.
\end{proof}

\begin{proposition}
\label{prop.case2.last}
Let $\gamma \in \widetilde{\mathcal{U}}^0$ satisfy either necessary condition of Proposition \ref{prop.case2.holo}.  Then $\gamma$ satisfies the sufficient conditions of Lemma \ref{fundconds} and $\langle \gamma \rangle$ acts freely, properly discontinuously, and cocompactly on $\Omega$.
\end{proposition}

\begin{proof}
First suppose that $\bar{\gamma} = q(\gamma)$ is a translation by a timelike vector $v$.  Considering $\gamma^{-1}$ if necessary, we may assume $v$ is future-pointing. 
Then evidently, for $q_0$ the origin of $\mbox{\bf Min}^-(\tilde{p}_0)$, we have $q_0\ll\gamma. q_0$. 

Under the Minkowski embedding $\iota$ in $\Ein$, the limit $\lim_{i \rightarrow \pm \infty} \iota(  iv )$ is $p_0$.  Thus $\lim_{i \rightarrow \pm \infty} \gamma^i.q_0 \in \{ \tilde{p}_i \}$.  By Corollary \ref{Minboundary}, both limits belong to $\{\alpha^{-1}.\tilde{p}_0, \tilde{p}_0, \alpha.\tilde{p}_0\}$.  On the other hand, the forward limit is in the future of $q_0$, so could only equal $\alpha.\tilde{p}_0$.  The backward limit is in the past of $q_0$.  Lemma \ref{Eincomplement} gives that the complement of $\mbox{\bf Min}^-(\tilde{p}_0)$ is $J^+(\tilde{p}_0) \cup J^-(\tilde{p}_0)$---that is, points of $\mbox{\bf Min}^-(\tilde{p}_0)$ are not causally related to $\tilde{p}_0$.   Then $\lim_{i \rightarrow - \infty} \gamma^i.q_0 = \alpha^{-1}.\tilde{p}_0$.  The sufficient condition (2) is thus verified for the case that $\bar{\gamma} \in \ker L$.

We now consider $\bar{\gamma} = L_{\bar{\gamma}} + v$, satisfying $L_{\bar{\gamma}} \in \mathcal{U}$ and $\langle  v, e_1 \rangle = v_n \neq 0$.  For $x \in \Min$, the linear part $L_{\bar{\gamma}} = U \in \mathcal{U}$ acts by
$$ U.x = x - x_n w + \left( \langle x, w \rangle - x_n ||w||^2/2 \right) e_1 \qquad \mbox{for some } w \in e_1^\perp \cap e_n^\perp$$
Replacing $\gamma$ by $\gamma^{-1}$ if necessary, we can arrange that $v_n > 0$.
In the Minkowski patch of $\Ein$, 
\begin{equation}
\label{eqn.case2.asymp}
 \bar{\gamma}^i.0 = \sum_{j=0}^{i-1} U^j.v = \sum_{j=0}^{i-1} \left( v - j v_n w + \left( j \langle v, w \rangle - j^2 v_n ||w||^2/2 \right) e_1\right) \qquad i \in \BN
\end{equation}

To evaluate the causal asymptotics of this sequence, recall $\partial_\theta$ given by (\ref{eqn.dtheta}).  Under the inverse of the Minkowski chart (see \S \ref{sec.geom.ein}), it pushes forward to 
$$ \left((\iota^{-1})_* \partial_\theta \right)_x = (-q_{n-1,1}(x)/2 - 1) (\partial_1 - \partial_n) + (x_n - x_1) (x_1 \partial_1 + \cdots + x_n \partial_n)$$
which has negative inner product with $\partial_n - \partial_1$ everywhere in $\Min$.  We may use this latter vector field for time orientation on $\Min$ (see \cite[Lem 5.29]{oneill.semi.riem}).  Since it is a constant vector field, a vector $x$ will be in the future light cone of the origin exactly when $\langle x, e_n - e_1 \rangle < 0$.

Now 
$$ \langle \bar{\gamma}^i.0, e_n - e_1 \rangle = - \frac{i(i-1)(2i-1) v_n ||w||^2}{12} + O(i^2) \rightarrow - \infty \mbox{ as } i \rightarrow \infty$$
Replacing $\gamma$ by $\gamma^k$ for some sufficiently large $k > 0$, we obtain $q_0 \ll \gamma.q_0$, and condition (1) of Lemma \ref{fundconds} is verified in this case.

The dominant term in (\ref{eqn.case2.asymp}) is the coefficient of $e_1$.  Under the Minkowski embedding, 
the limit points of $\{ \bar{\gamma}^i.\bar{q}_0 \ : \ i > 0 \}$ are in ${\bf P}(\mbox{span} \{e_0,e_1 \}) = \bar{\Delta}$ (here $\bar{q}_0 = \pi_{\mathrm{Ein}} (q_0) = \iota(0)$).  On the other hand, $\bar{\gamma}$ acts nontrivially on $\bar{\Delta}$ with unique attracting fixed point $p_0=[e_0]$.   Therefore $\lim_{i \rightarrow \infty} \bar{\gamma}^i.\bar{q}_0 = p_0$.  
For $i < 0$, the formula (\ref{eqn.case2.asymp}) becomes
$$\bar{\gamma}^i.0 = \sum_{j=1}^{-i} - U^{-j}.v = \sum_{j=1}^{-i} \left( - v - j v_n w + \left( j \langle v, w \rangle + j^2 v_n ||w||^2/2 \right) e_1\right) \qquad i < 0$$
Here again the $e_1$-component is dominant, and the same argument as for $i > 0$ gives $\lim_{i \rightarrow - \infty} \bar{\gamma}^i.\bar{q}_0 = p_0$.
  Then we conclude as in the translation case that $\lim_{i \rightarrow \pm \infty} \gamma^i.q_0 = \alpha^{\pm}.\tilde{p}_0$.  Sufficient condition (2) of Lemma \ref{fundconds} is verified.
\end{proof}

This case is now finished; the main results are summarized below.   The examples with $\gamma$ a timelike translation were previously discovered by C. Frances in his dissertation \cite[Sec 7.6.3]{frances.these}.

\begin{theorem}
\label{thm.case2}
If the image of $\bar{\delta}$ does not contain $p_0$ but does meet $\bar{\Delta}$, then, up to composition with a conformal transformation, $\delta$ is a diffeomorphism onto $\Omega = \mbox{\bf Min}^-(\tilde{p}_0) \cup L(\tilde{p}_0,\alpha.\tilde{p}_0) \cup \mbox{\bf Min}^+(\tilde{p}_0)$.  The holonomy image $\Gamma = \langle \gamma \rangle$ for $\gamma \in \widetilde{\mathcal{U}}^0$ satisfying the necessary and sufficient conditions in Proposition \ref{prop.case2.holo}.  The manifold $M$ is diffeomorphic to ${\bf S}^{n-1} \times {\bf S}^1$.
\end{theorem}

\section{The developing image does not meet the photon $\bar{\Delta}$ but does meet the lightcone $L(p_0)$.}
\label{sec.case3}

In this case, the
unipotence of the holonomy image leads, with the help of Proposition \ref{prop.fgh.technique},
to $\delta$ being a diffeomorphism of $\widetilde{M}$ onto
$\tEin\setminus\Delta$. 
Using, among other things, the algebraic hull of the holonomy, we prove that there are no examples in odd dimensions.  In even dimensions, we find a family of  Heisenberg fiber bundles over the circle.

\subsection{Development and holonomy for light cone components}
\label{sec.dev.lightcone.case3}

Let $\tilde{p}_0 \in \pi_{\mathrm{Ein}}^{-1}(p_0)$.  The space
$L(\tilde{p}_0)\setminus\Delta$ consists of infinitely
many connected components $\{ S_i \}$.  By case (2) of Proposition
\ref{prop.lightcone.complete}, $\delta$ maps each connected component
of $\delta^{-1}(S_i) \subset \widetilde{M}$ diffeomorphically to $S_i$ for each $i$.

We now choose components $S \subset L(\tilde{p}_0) \backslash
\Delta$ and $\Sigma \subset \delta^{-1}(S)$.  Denote
by $\hat{\Gamma}_\Sigma$ the stabilizer of $\Sigma$ in the group of deck transformations.  Because
$\left. \delta\right|_\Sigma$ is a diffeomorphism to its image, the holonomy
restricts to an isomorphism of $\hat{\Gamma}_\Sigma$ with its image, which will be denoted $\Gamma_\Sigma$.  The latter lies in the stabilizer of $S$, hence, under our assumptions, in $\widetilde{\mathcal{U}}^0$.  The restriction of
$\widetilde{\mathcal{U}}^0$ to $S$ is faithful, thus so is the restriction of
$\hat{\Gamma}_\Sigma$ to $\Sigma$.

\begin{lemma}
  \label{lem.sigma.cocompact}
The action of $\hat{\Gamma}_\Sigma$ on $\Sigma$ is cocompact.  The
cohomological dimension $\mbox{cd }\Gamma_\Sigma$ equals $n-1$.
\end{lemma}

\begin{proof}
The hypersurface $\Sigma$ is a connected component of the $\delta$-inverse image of the closed, $\Gamma$-invariant set $L(\tilde{p}_0)$.  
Then $\pi_M(\Sigma) = \Sigma / \hat{\Gamma}_\Sigma$ is closed by our standard arguments with Proposition \ref{prop.closed.invt}.
The second statement
follows from
$\Sigma \cong S \cong L(p_0) \backslash \bar{\Delta}
\cong \BR^{n-1}$.  
\end{proof}

Now we will focus on $\bar{\Gamma}_\Sigma =
q(\Gamma_\Sigma)$.  Note that $\pi_{\mathrm{Ein}} \circ \delta$ is
a diffeomorphism from $\Sigma$ to $\pi_{\mathrm{Ein}}(S) = L(p_0) \backslash
\bar{\Delta}$.  This means that $q\circ h$ is an isomorphism from
$\hat{\Gamma}_\Sigma$ to $\bar{\Gamma}_\Sigma$.

For $g \in \mathcal{U}$, recall the affine decomposition $g = L_g + u_g$ of \S \ref{sec.unip.subgp}, and that the projection
 $L(\mathcal{U}) < \SO(1,n-1)$ is abelian.  

    \begin{lemma}
      \label{lem.linears.span}
The linear projection $L(\bar{\Gamma}_\Sigma)$ spans the abelian Lie
group $L(\mathcal{U})$.  The image is discrete only if $\ker L
\cap \bar{\Gamma}_\Sigma$ is nontrivial.
\end{lemma}

\begin{proof}
   The leaf space for the foliation
  of $L(p_0) \backslash \{ p_0 \}$ by photons is identified with the round sphere
  ${\bf S}^{n-2}$, on which the $\mathcal{U}$-action factors through
  the projection $L$ and is conformal.  The quotient of $L(p_0)
  \backslash \bar{\Delta}$ by the photon foliation is thus identified
  with the punctured round sphere.  Under stereographic projection,
  this leaf space is conformal to $\Euc^{n-2}$.  The action of
  $L(\mathcal{U}) \cong \BR^{n-2}$ is by translations.  (This is the same action as appears in the proof of Proposition \ref{prop.lightcone.complete}.)

  By Lemma \ref{lem.sigma.cocompact}, $\hat{\Gamma}_\Sigma$ is cocompact on
  $\Sigma$.  Since $\Sigma$ is $(q \circ h)$-equivariantly
  diffeomorphic via $S$ to $L(p_0) \backslash \bar{\Delta}$, it
  follows that $\bar{\Gamma}_\Sigma$ acts cocompactly on $L(p_0)
  \backslash \bar{\Delta}$, and thus on the leaf space of the
  invariant foliation by photons.  A group of translations of $\Euc^{n-2}$ is
  cocompact if and only if it spans $\BR^{n-2}$.  Thus
  $L(\bar{\Gamma}_\Sigma)$ spans $L(\mathcal{U})$.
  
  If $\ker L \cap
  \bar{\Gamma}_\Sigma$ is trivial, then $\bar{\Gamma}_\Sigma$ maps
  isomorphically to its image in $L(\mathcal{U}) \cong \BR^{n-2}$.
  Thus $\bar{\Gamma}_\Sigma$ is a free abelian group in this case.  On
  the other hand, by Lemma \ref{lem.sigma.cocompact}, and because
  $\Gamma_\Sigma \cong \bar{\Gamma}_\Sigma$, it has
  cohomological dimension $n-1$.  Thus $\BZ^{n-1} \cong
  \bar{\Gamma}_\Sigma \cong L(\bar{\Gamma}_\Sigma) < \BR^{n-2}$.  It
  follows that $L(\bar{\Gamma}_\Sigma)$ is not discrete in this case.
  \end{proof}

  \begin{lemma}
    \label{lem.kernel.1d}
  The intersection
  $\ker L \cap \bar{\Gamma}_\Sigma$ is generated
by a single, possibly trivial, null translation in the center of
$\mathcal{U}$.  If it is nontrivial, then $L(\bar{\Gamma}_\Sigma)$ is discrete.
\end{lemma}

\begin{proof}
Denote $\mathcal{T} = \ker L \cap \bar{\Gamma}_\Sigma$.  It is a free
abelian group; let its rank be $k$.  The quotient $L(\bar{\Gamma}_\Sigma) <
L(\mathcal{U})$ is also a free abelian group, of rank at least $n-2$
by Lemma \ref{lem.linears.span}.  On the other hand, $\mbox{cd }
\bar{\Gamma}_\Sigma = \mbox{cd } \Gamma_\Sigma = n-1$. We conclude that
$k \leq 1$.

Of course $\mathcal{T}$ is normal in $\bar{\Gamma}_\Sigma$.  If
$\mathcal{T}$ is nontrivial, it is isomorphic to
$\BZ$. The conjugation action of $\bar{\Gamma}_\Sigma$ on $\mathcal{T}$ is
moreover unipotent; therefore, it is trivial, and $\mathcal{T}$ 
is central in $\bar{\Gamma}_\Sigma$.
Because $L(\bar{\Gamma}_\Sigma)$ spans $L(\mathcal{U})$,
$\mathcal{T}$ must be contained in the common
fixed space of $L(\mathcal{U})$ in $\BR^n$, which is a one-dimensional
isotropic subspace.
It acts on $\Ein$ by $\{ \tau^s \}$.

Finally, if $\mathcal{T} \neq 1$, then the rank of the free
abelian group $L(\bar{\Gamma}_\Sigma) \cong
\bar{\Gamma}_\Sigma/\mathcal{T}$ is at most $n-2$ by Lemma \ref{lem.sigma.cocompact}.  Again, because this
image spans $L(\mathcal{U})$, its rank is $n-2$ and it is discrete in
this case.
  \end{proof}



\begin{corollary} If $n=3$ or $\ker L \cap \bar{\Gamma}_\Sigma = 1$, then $\Gamma_\Sigma \cong \bar{\Gamma}_\Sigma$ is
  abelian.
\end{corollary}

\begin{proof}
If $\ker L \cap \bar{\Gamma}_\Sigma = 1$, then $
\bar{\Gamma}_\Sigma \cong L(\bar{\Gamma}_\Sigma)$, which is abelian.  If
$n=3$ and $\ker L \cap \bar{\Gamma}_\Sigma$ is nontrivial, then by
Lemma \ref{lem.kernel.1d}, $\bar{\Gamma}_\Sigma$ is a central extension of
$L(\bar{\Gamma}_\Sigma) \cong \BZ$ by $\BZ$,
which is necessarily abelian.  
\end{proof}

\subsection{Development and holonomy for Minkowski patches}
\label{sec.dev.min.case3}

  Let $\mathcal{H}$ be the foliation of $\tEin\backslash
\Delta$ by the fibers of
$\rho_\Delta$ (see \ref{sec.tau}).  Because the developing image is contained
in $\tEin \backslash \Delta$ and $\mathcal{H}$ is invariant by $\Gamma$, it pulls back to a foliation $\hat{\mathcal{H}}$ of
$\widetilde{M}$ by degenerate hypersurfaces, invariant by the group of deck transformations, which we will denote $\hat{\Gamma}$.
Each leaf of $\hat{\mathcal{H}}$ is a connected component of the
$\delta$-preimage of a leaf in $\mathcal{H}$.
One of these
leaves is $\Sigma$, which is already known to map
diffeomorphically under $\delta$ to $S$.  We now prove 
the same property for every leaf of the foliation.

\begin{proposition}
  \label{prop.leaf.diffeo}
For the foliation of $\widetilde{M}$ pulled back by $\delta$ from $\mathcal{H}$, each leaf in $\widetilde{M}$ maps diffeomorphically to a
leaf of $\mathcal{H}$.
\end{proposition}

\begin{proof}
We begin by defining vector fields $Y_2, \ldots, Y_{n-1}$ on $\Ein \backslash \bar{\Delta}$.  In the Minkowksi chart $\mbox{\bf Min}(p_0) = \mbox{\bf Min}([e_0])$, they are coordinate vector fields of 
$$ \iota : (y_1, y_2, \ldots, y_n) \mapsto [-\frac{q(y)}{2} : y_1: y_2: \cdots : y_n :1 ]$$
The two Minkowski charts $\mbox{\bf Min}([e_0])$ and $\mbox{\bf Min}([e_1])$ cover $\Ein \backslash \bar{\Delta}$---their complement is the projectivization of $e_0^\perp \cap e_1^\perp = \mbox{span} \{ e_0, e_1 \}$.
The second Minkowski chart is
$$ (x_1, \ldots, x_n ) \mapsto [x_1 : - \frac{q(x)}{2} : x_2 : \cdots : x_{n-1} : 1 : x_n]$$  
The change of coordinates is
$$ x_1 = - \frac{q(y)}{2 y_n} \qquad x_i = \frac{y_i}{y_n}, \ i=2, \ldots, n-1 \qquad x_n = \frac{1}{y_n}$$
For $i=2, \ldots, n-2$, the push forward by the change of coordinates is
$$ \frac{\partial}{\partial y_i } \mapsto x_n \frac{\partial}{\partial x_i } - x_i  \frac{\partial}{\partial x_1}$$
These vector fields thus extend to smooth vector fields $Y_2, \ldots, Y_{n-2}$ on $\Ein \backslash \bar{\Delta}$.

Note that $\mbox{\bf Min}([e_1]) \backslash \mbox{\bf Min}([e_0]) = L(p_0) \backslash \bar{\Delta}$ corresponds to $x_n=0$, and that here, the $Y_i$s are all tangent to the foliation $\bar{\mathcal{F}}$ of the lightcone by photons.  Then evidently in restriction to $\bar{S} = L(p_0) \backslash \bar{\Delta}$, the $\mathcal{U}$-action leaves $Y_i$ invariant modulo $\bar{\mathcal{F}}$ for all $i$.  Also, as they do not depend on $x_1$, the $Y_i$s are all projectable modulo $\bar{\mathcal{F}}$ on $\bar{S}$.
 In the coordinates on $\mbox{\bf Min}([e_0]) = \mbox{\bf Min}(p_0)$, the foliation $\bar{\mathcal{F}}$ by photons is the linear foliation tangent to  $e_1$, and $L(\mathcal{U})$ is trivial on $e_1^\perp/\BR e_1$, which means it fixes $Y_i = \frac{\partial}{\partial y^i}$ modulo $\bar{\mathcal{F}}$ for all $i$. As the $Y_i$ are all constant in these coordinates, they are projectable modulo $\bar{\mathcal{F}}$ inside $\mbox{\bf Min}(p_0)$, too.  Lastly, note that the $Y_i$s commute on $\Ein \backslash \Delta$.

Now lift $\{ Y_2, \ldots, Y_{n-1} \}$ to $\tEin \backslash \Delta$, keeping the same notation for them.  They are $\Gamma$-invariant modulo the foliation $\mathcal{F}$ by photons and projectable modulo $\mathcal{F}$.  
The $\tau$-flow is simply transitive on each photon of $\tEin \backslash \Delta$.  As always, since $\Gamma$ commutes with $\{ \tau^s \}$, it pulls back to a complete flow $\{ \hat{\tau}^s \}$ on $\widetilde{M}$ by Proposition \ref{prop.vf.pullback}.  Thus for the subfoliation $\hat{\mathcal{F}} \subset \hat{\mathcal{H}}$ by photons, each leaf maps diffeomorphically under $\delta$ to its image in $\tEin \backslash \Delta$.

As previously noted, every leaf of $\hat{\mathcal{H}}$ mapping into $L(\tilde{p}_0)$ maps diffeomorphically to a connected component of $L(p_0) \backslash \Delta$.  Therefore, we will consider $\mathcal{H}_y \subset \mbox{\bf Min}^-(\tilde{p}_j)$ for some $j$ for the remainder of the proof.

Part (1) of Proposition \ref{prop.fgh.technique} gives that the image of $\delta$ is invariant by the flow along any $Y \in \mbox{span} \{ Y_2, \ldots, Y_{n-1} \}$.  For $y \in \mbox{\bf Min}^-(\tilde{p}_j)$, the orbit of $y$ under these flows, saturated by $\mathcal{F}$, is the full leaf $\mathcal{H}_y$.  
For $y \in \mbox{\bf Min}^-(\tilde{p}_j)$, the $\{ Y_2, \ldots, Y_{n-1} \}$ form a framing modulo $\mathcal{F}$ in restriction to $\mathcal{H}_y$.  Then Proposition \ref{prop.fgh.technique} part (2), for $W = \mathcal{H}_y$, gives that the leaf $\hat{\mathcal{H}}_x$ mapping onto $\mathcal{H}_y$ maps by a covering.  As $\mathcal{H}_y \cong \BR^{n-1}$, this is a diffeomorphism.
  \end{proof}
  

\begin{proposition}
  \label{prop.min.diffeo}
Suppose that the image of $\delta$ intersects $\mbox{\bf Min}^-(\tilde{p}_i)$ for
some $\tilde{p}_i \in \pi_{\mathrm{Ein}}^{-1}(p_0)$.  Then an open
subset $\Omega_i \subset \widetilde{M}$ maps diffeomorphically under
$\delta$ to $\mbox{\bf Min}^-(\tilde{p}_i)$.
  \end{proposition}

\begin{proof}
By Proposition \ref{prop.leaf.diffeo}, the image of $\delta$ is $\mathcal{H}$-saturated; moreover, $\delta$ maps $\hat{\mathcal{H}}$-leaves diffeomorphically to $\mathcal{H}$-leaves.  To define a transverse vector field on $\tEin \backslash \Delta$, we reprise the Minkowski charts $\mbox{\bf Min}([e_0])$ and $\mbox{\bf Min}([e_1])$ on $\Ein$ as in the proof of that proposition.  

Under the change of coordinates, 
$$ \frac{\partial}{\partial y_n} \mapsto \left( \frac{q(x)}{2} + x_1 x_n \right) \frac{\partial}{\partial x_1} - \sum_{i=2}^{n-1} x_i x_n \frac{\partial}{\partial x_i} - x_n^2 \frac{\partial}{\partial x_n}.$$
Then the coordinate vector field $\partial/\partial y_n$ on $\mbox{\bf Min}(p_0)$ extends by $0$ to $L(p_0) \backslash \bar{\Delta}$, the set corresponding to $x_n = 0$ in $\mbox{\bf Min}([e_1])$.  We obtain a well-defined vector field $Y_n$ on $\Ein \backslash \bar{\Delta}$.  As $Y_n$ is constant in the chart $\mbox{\bf Min}(p_0)$ and $0$ elsewhere, it is projectable modulo $\bar{\mathcal{H}}$.  As $L(\mathcal{U})$ acts trivially on $\BR^{n-1,1}/e_1^\perp$, and $Y_n$ vanishes on the complement of $\mbox{\bf Min}(p_0)$, it is $\mathcal{U}$-invariant modulo $\bar{\mathcal{H}}$.

Now lift $Y_n$ to $\tEin \backslash \Delta$, where it will still be denoted $Y_n$.  It is $\Gamma$-invariant and projectable modulo $\mathcal{H}$.  Take $W = \mbox{\bf Min}^-(\tilde{p}_i)$. It is saturated by $\mathcal{H}$ and by the flow along $Y_n$; moreover, $Y_n$ is nonzero modulo $\mathcal{H}$ everywhere in $W$.  Then for $\Omega_i$ any connected component of $\delta^{-1}(W)$, Proposition \ref{prop.fgh.technique} says that $\Omega_i$ maps diffeomorphically onto $\mbox{\bf Min}^-(\tilde{p}_i)$.
  \end{proof}

A neighborhood of $\Sigma$ intersects at least one such $\Omega_i \subset \widetilde{M}$.  The holonomy subgroup $\Gamma_\Sigma$ lies in
$\widetilde{\mathcal{U}}^0$ and preserves each Minkowski patch in $\tEin$.
Thus $\hat{\Gamma}_\Sigma$ leaves invariant the open sets corresponding to Minkowski patches neighboring $\Sigma$ in $\widetilde{M}$.  Let
one such be $\Omega$.  It is mapped diffeomorphically under
$\bar{\delta}$ to $\mbox{\bf Min}(p_0)$.  The $\hat{\Gamma}_\Sigma$-action on $\Omega$
is thus conjugate to the $\bar{\Gamma}_\Sigma$-action on $\mbox{\bf Min}(p_0)$.
We now consider the latter action in detail, and recall that $\bar{\Gamma}_\Sigma \cong \Gamma_\Sigma$.  We will write the group $\Gamma_\Sigma$ and not distinguish from $\bar{\Gamma}_\Sigma$ for this discussion.

The leaf space of $\bar{\mathcal{H}}$ in $\mbox{\bf Min}(p_0)$ is
diffeomorphic via $\rho_{\bar{\Delta}}$ to $\bar{\Delta} \backslash \{ p_0 \} \cong \BR$.
Fix for the remainder of this section an identification $\mbox{\bf Min}(p_0)$ with $\Min$, with the quadratic form $q_{n-1,1}$ of \S \ref{sec.geom.ein}.
Denote by $D$ the
homomorphism sending an element of $\mathcal{U}$ to its action on this
leaf space.  The leaf space can also be identified with 
$\BR e_n$, on which $D(g)$ is translation by $\langle e_1, u_g \rangle$, 
the translational component of $g$ transverse to $\bar{\mathcal{H}}$.

\begin{proposition}
\label{prop.in.ker.D}
  If $n \geq 4$
  , then
  $\Gamma_\Sigma \subset \ker D$.  
On the leaf space of $\hat{\mathcal{H}}$ in $\widetilde{M}$,
  the action of $\hat{\Gamma}_\Sigma$ is trivial.
\end{proposition}

For $g,h \in \mathcal{U}$, the commutator $[g,h] = ghg^{-1}h^{-1}$ is easily seen to be
\begin{equation}
\label{eqn.commutator}
[ g,h]  =  (L_g - \mbox{Id})u_h - (L_h - \mbox{Id}) u_g
    \end{equation}

The identification of $L(\mathcal{U})$ with $\BR^{n-2}$ can be made
explicitly as follows:  
in our coordinates on $\Min$, elements $L_g \in L(\mathcal{U})$ act by
\begin{eqnarray}
\label{eqn.linear.action}
L_g(v) &  \equiv &  v - \langle v, e_1 \rangle \ell_g  \  \mbox{mod } \BR e_1 
\end{eqnarray}
    for a unique $\ell_g \in e_1^\perp/ \BR e_1  \cong \BR^{n-2}$.  Denote $V = e_1^\perp/\BR e_1$.  

    \begin{proof}
The aim is to prove that for all $g \in \Gamma_\Sigma$, the translational component  $u_g \in e_1^\perp$.  Since the linear part $L(g)$ preserves each leaf of $\bar{\mathcal{H}}$ in $\mbox{\bf Min}(p_0)$, it will then follow that $g$ preserves each leaf.  The corresponding statement about $\hat{\Gamma}_\Sigma$ will follow from Proposition \ref{prop.min.diffeo}.
    
Let $g,h \in \Gamma_\Sigma$, and let $\ell_g, \ell_h \in
V$ correspond to $L(g)$ and $L(h)$ as above.  The commutator
given by (\ref{eqn.commutator}) belongs to
$\ker L \cap \Gamma_\Sigma$.
The latter subgroup acts on $\Min$ by a group of translations in the
direction of $e_1$ by Lemma \ref{lem.kernel.1d}.  Thus
$$ (L_g - \mbox{Id})u_h - (L_h - \mbox{Id}) u_g \equiv \langle u_h,
e_1 \rangle \ell_g - \langle u_g, e_1 \rangle \ell_h \equiv 0 \qquad \mbox{mod
} \BR e_1$$
for all $g,h \in
\Gamma_\Sigma$.  If for one $g \in \Gamma_\Sigma$, the
component $\langle u_g , e_1 \rangle \neq 0$, then $\ell_g \neq
0$ by Lemma \ref{lem.kernel.1d}, and $\ell_h$ must be
linearly dependent with $\ell_g$ for all $h \in \Gamma_\Sigma$.
If $n \geq 4$, then the latter conclusion contradicts Lemma \ref{lem.linears.span}.
\end{proof}

The action of $\Gamma_\Sigma$ on $\Min$ is affine and properly discontinous. 
Viewing $\mathcal{U}$ as a subgroup of the affine group of $\Min$, there is a connected algebraic hull $H < \mathcal{U}$ which acts properly on $\Min$ and contains $\Gamma_\Sigma$ as a cocompact lattice (see \cite[Thm 1.4]{fried.goldman.3d}).  In fact, $H$ acts freely on $\Min$ as well \cite[Lem 1.9]{fried.goldman.3d}.  By the proposition above, $H$ is contained in $\ker D$ if $n \geq 4$, so it acts freely and properly on each hyperplane in the foliation $\bar{\mathcal{H}}$ in this case.

\begin{lemma}
\label{lem.H.kernel}
For $H$ the algebraic hull of $\Gamma_\Sigma$ in $\mbox{Aff}(\Min)$, there are nontrivial translations in $H$, that is, $\ker L \cap H \neq 1$. 
\end{lemma}

\begin{proof}
Suppose that $\ker L \cap \Gamma_\Sigma = 1$.  By Lemma \ref{lem.sigma.cocompact} and because $L(\mathcal{U}) \cong \BR^{n-2}$, the group $\Gamma_\Sigma \cong \BZ^{n-1}$.
Thus the algebraic hull $H \cong \BR^{n-1}$, and it intersects the kernel of $L$ nontrivially.

If $\ker L \cap \Gamma_\Sigma \neq 1$, then evidently $H$ also has nontrivial intersection. 
\end{proof}

We establish some basic structural features of $H$.  Since $L(\Gamma_\Sigma)$ spans $L(\mathcal{U})$ by Lemma \ref{lem.linears.span}, the linear projection of the algebraic hull $L(H)$ equals $L(\mathcal{U})$.  For $h \notin \ker L$, the translational component $u_h$ is given by a 1-cocycle $L(H) \rightarrow \BR^{n-1,1} = \ker L$; by Proposition \ref{prop.in.ker.D}, this cocycle has values in $e_1^\perp$.  Because $L(H)$ acts trivially on $V = e_1^\perp / \BR e_1$, composing $L_h \mapsto u_h$ with projection to $V$ yields a homomorphism.
Recall the identification of $L(\mathcal{U})$ with $V$ via $\ell_g$ from (\ref{eqn.linear.action}).
Define 
$$\theta \in \mbox{End }V \qquad  \theta : \ell_h \mapsto \bar{u}_h  = [u_h] \ \mbox{mod } \BR e_1$$


\begin{proposition}
\label{prop.real.eigenvalues}
Assume that $\Gamma_\Sigma < \ker D$ and let $H$ be the algebraic hull, with associated $\theta \in \mbox{End }V$ as above.
 If $\theta$ has a real eigenvalue, then $H$ does not act properly on $\Min$; more precisely, if $r$ is a real eigenvalue, then $H$ acts with noncompact stabilizers on the affine hyperplane defined by $\langle x, e_1 \rangle = r$.
\end{proposition}

\begin{proof}
Let $\ell \in V$ be an eigenvector of $\theta$ with eigenvalue $r$.  Because $H$ projects onto $L(\mathcal{U})$, there is an element $h \in H$ with $\ell_h = \ell$.  By definition, $\theta(\ell_h) = r \ell_h \equiv [u_h] \ \mbox{mod } \BR e_1$.  Given $x \in \Min$, equation (\ref{eqn.linear.action}) says
$$ h(x) = L_h(x) + u_h \equiv x - \langle x, e_1 \rangle \ell_h + u_h \ \mbox{mod } \BR e_1$$
Thus $h$, and the 1-parameter subgroup of $H$ containing $h$, act trivially modulo $\BR e_1$ on the affine hyperplane defined by  $\langle x, e_1 \rangle = r$ .  The translations $\BR e_1$ belong to $H$, as well, by Lemma \ref{lem.H.kernel}.  Therefore any vector in this hyperplane has a noncompact stabilizer, proving nonproperness.
\end{proof}

\subsection{The global picture}
As above, $H$ is the algebraic hull of $\Gamma_\Sigma$ in $\mathcal{U}$.

\begin{proposition}
\label{prop.n=3}
  Suppose that $n=3$. 
  Then $\hat{\Gamma}_\Sigma$ does not act properly on $\widetilde{M}$. 
\end{proposition}

From this proposition we conclude there are no $3$-dimensional examples for which $\mbox{im } \bar{\delta}$ meets $L(p_0)$ but not $\bar{\Delta}$.

\begin{proof}
Let $S$ be a connected component of $L(p_0) \backslash \Delta$ in the image of $\delta$, as above.  Let $\mbox{\bf Min}^-(\tilde{p}_i)$ have $S$ in its closure, so that by Proposition \ref{prop.min.diffeo}, there is $\Omega_i$ mapping diffeomorphically to $\mbox{\bf Min}^-(\tilde{p}_i)$.

Let $\ell = \BR e_3 \subset \mbox{\bf Min}^{2,1}$.  Note that, under the Minkowski chart $\mbox{\bf Min}^{2,1} \rightarrow \mbox{\bf Min}(p_0)$ the closure $\overline{\ell}$ in $\mbox{\bf Ein}^{2,1}$ is a photon not meeting $\bar{\Delta}$.  
Now consider the lift $\ell \subset \mbox{\bf Min}^-(\tilde{p}_i)$; the closure $\overline{\ell}$  in $\widetilde{\mbox{\bf Ein}}^{2,1} \backslash \Delta$ is compact, contained in $\mbox{\bf Min}^-(\tilde{p}_i) \cup S$.  It is therefore the diffeomorphic image of a compact photon $C$ in $\widetilde{M}$.


Fix a generater $\ell$ of $e_1^\perp/\BR e_1$ and write $\ell_g = s_g \ell$ with $s_g \in \BR$ for each $g \in H$.  We can write $u_g \equiv a_g \ell + b_g e_3$ mod $\BR e_1$.  Then, whenever $s_g \neq 0$,  we can take  $t = a_g/s_g$ to obtain
$$ g(t e_3) = L_g(te_3) + u_g \equiv t e_3 - ts_g \ell + a_g \ell + b_g e_3 \equiv (t + b_g) e_3 \ \mbox{mod } \BR e_1$$

Because $H$ projects onto $L(\mathcal{U})$, there is an unbounded set of $h \in H$ with $s_h \neq 0$.  On the other hand, by Lemmas \ref{lem.H.kernel} and \ref{lem.kernel.1d}, the translations by $\BR e_1$ are in $H$.  Therefore, there is an unbounded subset of $h \in H$ for which $h.\bar{\ell} \cap \bar{\ell} \neq \emptyset$.  As $\Gamma_\Sigma$ is a lattice in $H$, there are infinitely many $\gamma \in \Gamma_\Sigma$ and a compact subset $D \subset H$ such that $\gamma.K \cap K \neq \emptyset$ for $K = D.\overline{\ell} \subset \mbox{\bf Min}^-(\tilde{p}_i) \cup S$.  Then there are a compact subset $\hat{K} \subset \widetilde{M}$ and infinitely many $\hat{\gamma} \in \hat{\Gamma}_\Sigma$ such that $\hat{\gamma}.\hat{K} \cap \hat{K} \neq \emptyset$.
\end{proof}

\begin{proposition}
  \label{prop.dev.case3}
  The developing map $\delta$ is a diffeomorphism of $\widetilde{M}$
  with $\tEin\backslash \Delta$.
\end{proposition}

This lemma will be used in the proof.

\begin{lemma}
  \label{lem.homeo.leaf.space}
The developing map $\delta$ induces a local homeomorphism between the
leaf space $(\tEin\backslash \Delta)/ \mathcal{H} \cong \Delta$ and the
corresponding leaf space in $\widetilde{M}$. In fact, the induced map
on leaf spaces
is a diffeomorphism to its image, which is an open interval in $\Delta$.
  \end{lemma}

  \begin{proof}
Let $p \in \widetilde{M}$ and let $W$ be a neighborhood of $p$ mapping
diffeomorphically under $\delta$ to its image in $\tEin$.  Shrink $W$ as necessary so
that $\delta(W)$ admits a transversal submanifold $T$ to $\mathcal{H}$
through $\delta(p)$ mapping diffeomorphically to its image $\overline{T}$ in the leaf
space of $\tEin\backslash \Delta$.  

Let $\tau$ be the corresponding transversal to the foliation in $W$,
so that $\delta$ restricted to $\tau$ is a diffeomorphism to $T$.  Let
$\bar{\tau}$ be the image of $\tau$ in the leaf space of
$\widetilde{M}$.  The map sending leaves of $\bar{\tau}$ to leaves of
$\bar{T}$ is well-defined because $\delta$ intertwines the two foliations.
The inverse map is the composition $\bar{T} \rightarrow T \rightarrow \tau \rightarrow \bar{\tau}$
of two differomorphisms with a quotient map; it is evidently continuous.

It remains to verify that the map $\bar{\tau} \rightarrow
\bar{T}$ is continuous.  An open subset $\bar{U}
\subset \bar{T}$ is diffeomorphic to a relatively open subset $U
\subset T$.  There is moreover an $\mathcal{H}$-saturated open subset
$\tilde{U}$ such that $\tilde{U} \cap T = U$.  Now
$$(\left. \delta \right|_W)^{-1}(U) = (\left. \delta
\right|_W)^{-1}(\tilde{U} \cap T) = \delta^{-1}(\tilde{U}) \cap \tau$$
The projection of the set on the left-hand side is the inverse image
of $\bar{U}$ in $\bar{\tau}$.  The expression on the
right-hand side exhibits it as an open subset
of $\bar{\tau}$.

Now that we have shown that the map on leaf spaces is a local
homeomorphism, it follows that the leaf space in $\widetilde{M}$ is a
one-dimensional manifold, and the map between leaf spaces is a local diffeomorphism.  The image in the leaf space of
$\tEin\backslash \Delta$ is open and connected, thus an
open interval.  The leaf space in $\widetilde{M}$ is therefore
diffeomorphic to $\BR$ and maps diffeomorphically to its image.
    \end{proof}

    We now prove Proposition \ref{prop.dev.case3}.
    
\begin{proof}
The image of the developing map $\delta$ is an open, contiguous union of
Minkowski charts $\mbox{\bf Min}(\tilde{p}_i)$ and their interstitial light
cone components $S_j$.  
By Lemma \ref{lem.homeo.leaf.space}, $\delta$ factors through a map on
leaf spaces which is a diffeomorphism to an open interval $I \subset
\Delta \cong \BR$.  The fibers of the composition $\widetilde{M}
\rightarrow I$ are the leaves of $\hat{\mathcal{H}}$,
which is thus a simple foliation of $\widetilde{M}$.  By Proposition \ref{prop.leaf.diffeo}, it
follows that $\widetilde{M}$ is diffeomorphic to $\BR^n$ and maps
diffeomorphically under $\delta$ to its image.  Therefore
$\hat{\Gamma}$ maps isomorphically under $h$ to its image $\Gamma < \widetilde{\mathcal{U}}$.

It remains to show
the image of $\delta$ is all of $\tEin \backslash \Delta$.
Suppose
that $I$ is a finite interval in $\BR$.  Then it contains $S_i$ for
finitely many $i$.  A finite index subgroup of $\Gamma$
stabilizes each light cone component; let one of them be $S$.  For
$\Sigma = \delta^{-1}(S)$, the group $\hat{\Gamma}_\Sigma$ is isomorphic
to $\mbox{Stab}(S)$ and has finite index
in $\hat{\Gamma}$.  Then $\mbox{cd } \Gamma = \mbox{cd } \Gamma_\Sigma =
n-1$ by Lemma \ref{lem.sigma.cocompact}.  This contradicts $\hat{\Gamma}$ acting properly
discontinuously and cocompactly on $\widetilde{M} \cong \BR^n$.  Therefore $I$ is infinite.  

Denote by $\widetilde{D} : \widetilde{\mathcal{U}} \rightarrow \widetilde{\mbox{PSL}}(2,\BR)$ the homomorphism given by restriction to $\Delta$, which also corresponds to the action on the leaf space of $\mathcal{H}$.  As $\widetilde{D}(\Gamma)$ acts cocompactly on $I$, it contains elements with nontrivial translation number.  The invariance of $I$ by such an element implies that $I = \BR$.
\end{proof}

By Proposition \ref{prop.n=3} we may assume $n \geq 4$, and then Proposition \ref{prop.in.ker.D} says that $\Gamma_\Sigma = \ker \widetilde{D}$, where $\widetilde{D}$ is as defined in the proof above.  The image $\widetilde{D}(\Gamma)$ projects under $\widetilde{\mbox{PSL}}(2,\BR) \rightarrow \mbox{PSL}(2,\BR)$ to the unipotent subgroup fixing $p_0$, which we will denote $\mathcal{V}$.  As with $\mathcal{U}$, this group has a lift to $\widetilde{\mbox{PSL}}(2,\BR)$, so $\widetilde{D}(\Gamma) < \widetilde{\mathcal{V}} \cong \widetilde{\mathcal{V}}^0 \times \BZ$.  But $\widetilde{D}(\Gamma) \cong \Gamma /\Gamma_\Sigma$ cannot intersect $ \widetilde{\mathcal{V}}^0 < \mbox{Stab}(p_0)$: the inverse image in $\Gamma$ of this intersection would stabilize $S$ and would thus be contained in $\Gamma_\Sigma$, a contradiction.  Therefore $\widetilde{D}(\Gamma) \cong \Gamma / \Gamma_\Sigma$ is generated by a single element 
acting properly discontinuously and cocompactly on the leaf space of $\hat{\mathcal{H}}$.

\begin{proposition}
The holonomy image $\Gamma$ is generated by $\Gamma_\Sigma$ and another element of the form $\alpha^{i_D} g_D$, where $i_D \in \BZ$ and $g_D \in \widetilde{\mathcal{U}}^0$ normalizes $\Gamma_\Sigma$.    
The projection $\bar{g}_D$ of $g_D$ to $\mathcal{U}$ belongs to $\ker D$.  
\end{proposition}

\begin{proof}
It remains only to verify the last claim.
Equation (\ref{eqn.commutator}) gives for $h \in \bar{\Gamma}_\Sigma$,
$$ \bar{g}_D h \bar{g}_D^{-1} = [\bar{g}_D,h]h = L_h + L_D (u_h) - (L_h - \mbox{Id}) u_D $$
The result is in $\bar{\Gamma}_\Sigma$ with image under $L$ equal $L_h$.  The translational part is
$$ L_D(u_h) - (L_h - \mbox{Id}) u_D \equiv u_h +  \langle u_D, e_1 \rangle \ell_h \ \mbox{mod } \BR e_1$$
which must equal $\bar{u}_h$ because $\ker L \cap \bar{\Gamma}_\Sigma < \BR e_1$.  Therefore $\langle u_D, e_1 \rangle = 0$ and $g_D \in \ker D$.

\end{proof}


%

We summarize the results obtained thus far, under the standing assumptions of this section: The developing map is a diffeomorphism to $\tEin\backslash \Delta$. The dimension of $M$ is $n\geq 4$.  The group $\hat{\Gamma}$ has a normal subgroup $\hat{\Gamma}_\Sigma$ which maps isomorphically under $q \circ h$ to its image $\bar{\Gamma}_\Sigma$, which has algebraic hull contained in $\ker D$.  Any element of $\Gamma \backslash \Gamma_\Sigma$ is of the form $\alpha^i g$ with $\bar{g} = \pi_{\mathrm{Ein}}(g) \in \ker D$ and $i \neq 0$.  The algebraic hull $H$ of $\bar{\Gamma}_\Sigma$ contains the center of $\mathcal{U}$ and is encoded by a linear endomorphism $\theta \in \mbox{End }V$
with no real eigenvalues; in particular $n - 2 = 2k$ for some $k \in \BN$.

\subsection{Heisenberg examples and conclusion}

The first result of this section underlies the construction of actions whenever the conditions summarized at the end of the previous section are fulfilled.  Then we give the classification for this case.

\begin{proposition}
\label{prop.heis.examples}
Let $n = 2(k+1)$.
Let $z_1, \bar{z_1}, \ldots, z_k, \bar{z_k} \in \BC \backslash \BR$ and let $\theta \in \GL(2k,\BR)$ be any element with these eigenvalues.  Define $H < \Aff \Min$ to be the connected group generated by
$$Z(\mathcal{U}) \cup \{ L_h + u_h :  u_h = \theta(\ell_h) \in \BR^{n-2}  \cong e_1^\perp \cap e_n^\perp , \ L_h \in L(\mathcal{U}) \}$$
Then $H$ acts simply transitively on $L(p_0) \backslash \bar{\Delta}$ and on $\mathcal{H}_r = \{ x \in \Min : \langle x, e_1 \rangle = r \}$ for each $r \in \BR$.
 \end{proposition}
 
 
\begin{proof}
Fix $r \in \BR$ and let $x + r e_n \in \mathcal{H}_r$ where $x \in e_1^\perp$.  An element $h \in H$ acts by
$$ h(x + r e_n) \equiv r e_n + x + u_h - r \ell_h \ \mbox{mod } \BR e_1$$
The result is congruent to $x + r e_n$ only if $u_h \equiv r \ell_h$.  By the construction of $H$, the latter condition implies $u_h \equiv \ell_h \equiv 0$ and $h \in Z(\mathcal{U})$.  Thus the stabilizer in $H$ of $x + re_n$ is trivial.  Moreover, for each $r \in \BR$, the set of $\bar{u}_h - r \ell_h = (\theta -r ) \ell_h$ spans $V$.  The action is simply transitive on $\mathcal{H}_r$.  

Recall the parametrization of $L(p_0) \backslash \bar{\Delta}$  by $\BR \times \BR^{n-2}$ of (\ref{eqn.param.lightcone}).
Assuming $u_h \in e_1^\perp$, the action of $L_h + u_h$ maps $(t,x)$ to 
$$  (t + \langle \bar{u}_h, x \rangle + \langle u_h, e_n \rangle,  x - \ell_h  )$$
The stablizer of $(t,x)$ has $\ell_h = 0 = \bar{u}_h$ which means $h \in Z(\mathcal{U})$.  The stabilizer also has $\langle u_h, e_n \rangle = 0$, which means $h$ is trivial.  It is also evident from the above formula that $H$ acts transitively on $\BR \times \BR^{n-2} \cong L(p_0) \backslash \bar{\Delta}$.
\end{proof}

From such an $H$-action, we can easily construct geometries on compact $M^{2k+2}$ satisfying the assumptions of this section.

\begin{theorem}
If the image of $\bar{\delta}$ does not meet $\bar{\Delta}$ but does meet $L(p_0)$, then $n$ is even, and $\delta$ is a diffeomorphism onto $\tEin \backslash \Delta$.  For a nilpotent group $H$ as in Proposition \ref{prop.heis.examples} and a lattice $\Gamma' < H < \widetilde{\mathcal{U}}^0$, the holonomy is a nilpotent extension $\Gamma = \langle \alpha^i g \rangle \ltimes \Gamma'$, for some $i > 0$ and $g \in \ker (D \circ q) < \widetilde{\mathcal{U}}^0$.
In this case, $M = (\tEin \backslash \Delta) / \Gamma$ is a nilmanifold of degree at most 3, which fibers over $S^1$ with degree-2 nilmanifold fibers.
\end{theorem}

\begin{proof}
Let $\Gamma_\Sigma$ be as defined in \S \ref{sec.dev.lightcone.case3}, and let $H$ be its algebraic hull in $\mathcal{U}$.  By Proposition \ref{prop.n=3} $n \neq 3$.  By Proposition \ref{prop.in.ker.D} and Lemmas \ref{lem.kernel.1d} and \ref{lem.H.kernel} , $H < \ker D$ and $Z(\mathcal{U}) < H$, so $H$ is determined by an endomorphism $\theta$ of $\BR^{n-2}$ as in \S \ref{sec.dev.min.case3}.  By Proposition \ref{prop.real.eigenvalues}, $\theta$ has no real eigenvalues, which forces $n-2$, hence $n$, to be even.  We conclude that $H$ must be as in Proposition \ref{prop.heis.examples}.  It is of nilpotence degree 2.
 
Proposition \ref{prop.dev.case3} establishes that $\delta$ is a diffeomorphism onto $\tEin \backslash \Delta$.

The conclusions of Proposition \ref{prop.heis.examples} allow us to identify the fibration $\rho_\Delta : \Ein \backslash \Delta \rightarrow \Delta$ with a principal $H$-bundle.  The image of $H$ under the isomorphism $\mathcal{U} \rightarrow \widetilde{\mathcal{U}}^0$ will also be denoted $H$.  Then $\tEin \backslash \Delta \cong \widetilde{M}$ is a principal $H$-bundle over $\Delta \cong \BR$.  

By construction, $\Gamma_\Sigma$ is a cocompact lattice in $H$.  Let  $\Gamma'$ be the image of $\Gamma_\Sigma$ in $\widetilde{\mathcal{U}}^0$.  As observed after Proposition \ref{prop.dev.case3}, $\Gamma/\Gamma'$ is generated by an element of the form $\alpha^i g$ with $i \neq 0$ and $g$ in the normalizer in $\widetilde{\mathcal{U}}^0$ of $\Gamma'$; we may assume $i > 0$. Since $\widetilde{\mathcal{U}}$ splits as a product $\BZ \times \widetilde{\mathcal{U}}^0$, there is a splitting $\Gamma = \langle \alpha^i g \rangle \ltimes \Gamma'$.  The remaining conclusions follow. 
\end{proof}


\section{Conclusion of classification, including essential examples}
\label{sec.case4.essential}

The last case in our classfication reduces to a class of manifolds that have been well studied.

\begin{theorem}
If the image of $\bar{\delta}$ does not meet $L(p_0)$, then $M$ is a complete $(\OO(n-1,1) \ltimes  \BR^{n-1,1}, \Min)$-manifold.  It is $N/\Gamma$ for $N <  \OO(n-1,1) \ltimes \BR^{n-1,1}$ a nilpotent group of degree at most 3 acting simply transitively on $\Min$, and $\Gamma$ a cocompact lattice in $N$.
\end{theorem}

\begin{proof}
The hypothesis implies that $\bar{\delta}$ maps into $\mbox{\bf Min}(p_0)$.  The $\mathcal{U}$-action on here is by affine isometries.  As $q(\Gamma) < \mathcal{U}$, the pair $(\bar{\delta}, q \circ h)$ defines an $(\OO(n-1,1) \ltimes  \BR^{n-1,1}, \Min)$-structure on $M$.  Y. Carri\`ere proved that any such structure on a closed manifold is complete---that is, $\bar{\delta}$ is a diffeomorphism of $\widetilde{M}$ onto $\Min$ \cite{carriere.lorentz.complete}.  

Now $h$ is an isomorphism to $\Gamma$, and $M \cong \Min/\Gamma$.  The algebraic hull of $\Gamma$ is a unipotent subgroup of $\OO(n-1,1) \ltimes  \BR^{n-1,1}$ acting simply transitively on $\Min$.  These were classified by F. Gr\"unewald and G. Margulis in \cite[Thm 1.8]{gm.min.unip}.  They are all nilpotent groups of degree at most $3$.
\end{proof}

\subsection{Proof of Theorem \ref{thm.essential}.}

Recall that the flow $\{ \tau^s \}$ is in the center of $\mathcal{U}$.  That means it always descends to a conformal flow on $M$, which we will denote $\{ \bar{\tau}^s \}$, in our setting of unipotent holonomy.  

\begin{proposition}
\label{prop.tau.essential}
For $M$ as in cases (1) or (2) of Theorem \ref{mainthm}, the conformal flow $\{ \bar{\tau}^s \}$ is essential.  For $M$ as in cases (3) or (4), it is inessential.
\end{proposition}

\begin{proof}
In cases (1) or (2), let $I \subset \Delta$ be a nontrivial, open interval contained in the image of $\delta$.  The inverse image $\rho_\Delta^{-1}(I) \subset \tEin \backslash \Delta$ is open.  For any volume $\nu$ on $\delta(\widetilde{M})$, the volume $\nu(I) = 0$ while $\nu(\rho_\Delta^{-1}(I) \cap \delta(\widetilde{M})) \neq 0$.  In fact, for any compact set $K \subset \rho_\Delta^{-1}(I)$ with nonempty interior, $\nu(K) \neq 0$.  Now $\tau^t(K)$ tends uniformly to a subset of $I$ ast $t \rightarrow \infty$, which implies that $\tau$ does not preserve any volume on $\delta(M)$.  In particular, it does not preserve a $\Gamma$-invariant volume lifted from $M$.  Therefore $\{ \bar{\tau}^s \}$ does not preserve any volume on $M$, so it must be essential.

In case (3), $H$ is not abelian.  The lattice $\Gamma' = \Gamma_\Sigma$ having $H$ as algebraic hull is therefore not abelian, so its commutator subgroup intersects $Z(\mathcal{U}) = \{ \tau^s \}$ nontrivially.  The $\{ \bar{\tau}^s \}$-flow on $M$ factors through the quotient by this intersection, which is ${\bf S}^1$.  It is therefore not essential.

Finally, $\{ \tau^s \} < \mathcal{U}$ acts on $\Min$ by a lightlike translation, which is isometric.  In case (4), $\Gamma < \mathcal{U}$, and $\mathcal{U}$ is isometric on $\Min$, so the flow $\{ \bar{\tau}^s \}$ is isometric on $M$.
\end{proof}

In light of Theorem \ref{mainthm} and Proposition \ref{prop.tau.essential}, it remains only to prove that $M$ as in cases (3) or (4) does not admit an essential conformal flow.  

In case (4), any conformal flow on $M$ lifts on $\widetilde{M} \cong \Min$ to $\Conf \Min \cong \CO(n-1,1) \ltimes \BR^{n-1,1}$.  Any non-isometric flow in this group is a homothety.  But a homothety cannot descend to a closed manifold.  Thus there is no essential conformal flow in this case.

In case (3),
any conformal flow of $\delta(\widetilde{M})$ that descends to $M$ belongs to the identity component of the normalizer of $\Gamma$ and to the stabilizer of $\Delta$.  Because $\Gamma$ is discrete, such a flow belongs in fact to $C^0(\Gamma)$, the identity component of the centralizer of $\Gamma$. 
Now $\Gamma = \langle \alpha^i g \rangle \ltimes \Gamma'$, so $C^0(\Gamma) < C^0(\Gamma') = C^0(H)$, as $H$ is the Zariski closure of $\Gamma$.  We descend to $\Ein$ and analyze the centralizer of $H$ intersect the stabilizer of $\bar{\Delta}$ in $G$.  Let $n=2k+2$.

The stabilizer of $\bar{\Delta}$ in $G$ corresponds to the stabilizer $P$ of $\mbox{span} \{ e_0,e_1 \}$, and has Levi decomposition 
$$P \cong (GL(2,\BR) \times O(n)) \ltimes \mbox{Heis}_{2n-3}$$
The kernel of $D$ intersect $\mathcal{U}$ is precisely $\mbox{Heis}_{2n-3}$.  
It contains $H$, and their centers coincide; we will denote this central subgroup by $Z$.  
The usual generators of $\mbox{Heis}_{2n-3}$ are $\{ x_j, y_j, \tau :  [x_j,y_j]=\tau, \ j=1, \ldots, n-2 \}$.  These can be chosen with $y_j \in \ker L$ for all $j$, so that we can view $\{x_j \}$ as a basis for $L(\mathcal{U})$ and $\{ y_j \}$ as a basis for $V = e_1^\perp/\BR e_1$.  The isomorphism $L(\mathcal{U}) \cong V$ of (\ref{eqn.linear.action}) corresponds in this basis to $x_j \mapsto y_j$ for all $j$.

Since $C^0(H)$ must act trivially on $Z$, its projection to $\GL(2,\BR)$ has image in $\SL(2,\BR)$.  Now we consider the centralizer in $\SL(2,\BR)$ of $H/Z < \mbox{Heis}_{2n-3}/Z \cong \BR^{2(n-2)} = \BR^{4k}$.  The latter representation is $\oplus_{j=1}^{2k} E_j$, where each $E_j = \mbox{span} \{ x_j,y_j \}$ is a copy of the standard representation.  
Recall that $H/Z$ is the graph of an isomorphism $\theta : \mbox{span} \{ x_j \} \rightarrow \mbox{span} \{ y_j \}$ with no real eigenvalues, as in Proposition \ref{prop.heis.examples}.  If $\{ g^t \} < \SL(2,\BR)$ is a noncompact 1-parameter subgroup centralizing $H$, it leaves invariant the graph of $\theta$, which we will denote $W$.  Since $W \cap \mbox{span} \{ y_j \} = 0$, each $W \cap E_j$ is of dimension at most 1.  Let $v^\pm$ denote the eigenvectors of $\{ g^t \}$ on $\BR^2$, possibly equal; denote by $v_j^\pm$ the corresponding vectors in $E_j$ for each $j$. Invariance implies that $W$ equals the sum of its intersections with $\BR v_j^\pm$.   
Any nontrivial intersection is generated by $x_j + \alpha y_j$ for some $\alpha \in \BR$, which would correspond to a real eigenvalue $\alpha$ of $\theta$, a contradiction.  We conclude that the centralizer of $H$ in $\SL(2,\BR)$ is compact.

Now we consider the intersection of $C^0(H)$ with the unipotent radical $\mbox{Heis}_{2k-3} \cong \ker D \cap \mathcal{U}$.   It follows from $H$ projecting onto $L(\mathcal{U})$ (Lemma \ref{lem.linears.span}) that this centralizer equals $Z$, which is $\{ \tau^s \}$.  By Proposition \ref{prop.tau.essential}, there is no essential conformal flow in this case, and the proof is complete.

\bibliographystyle{amsplain}
\bibliography{karinsrefs}

 \begin{tabular}{lll}
 Rachel Lee & \qquad & Karin Melnick   \\
  Department of Mathematics  &\qquad & Department of Mathematics \\
  4176 Campus Drive & \qquad & 6 avenue de la Fonte \\
 University of Maryland & \qquad & University of Luxembourg \\
College Park, MD 20742  & \qquad &  L-4364 Esch-sur-Alzette  \\
  USA &  \qquad & Luxembourg \\
rachel.youcis.lee@gmail.com &  \qquad & karin.melnick@uni.lu
 \end{tabular}

\end{document}